\documentclass{article}

 \PassOptionsToPackage{numbers, compress}{natbib}

 \usepackage[preprint]{neurips_2019}

\usepackage{lscape} 

\usepackage{graphicx}
\usepackage{subfigure}

\usepackage{graphicx} 
\usepackage{wrapfig} 
\usepackage{float}	  
\usepackage{calc} %
\usepackage{caption}
\usepackage{natbib} 

\usepackage{fontawesome} 
\usepackage[nolist,smaller]{acronym}

\usepackage{algorithm}
\usepackage[noend]{algpseudocode}

\usepackage{amsfonts, amsmath, amsthm} 
\usepackage{cancel} 
\usepackage{bbm, bm} 
\usepackage{IEEEtrantools}
\theoremstyle{definition}
\newtheorem{theorem}{Theorem} 
 
\newtheorem{assumption}{Assumption} 
\newtheorem{definition}{Definition}
\newtheorem{corollary}{Corollary}
\newtheorem{proposition}{Proposition}
\newtheorem{remark}{Remark}
\newtheorem{lemma}{Lemma}
\usepackage[dvipsnames]{xcolor}
\usepackage{enumitem}

\definecolor{navy}{RGB}{0, 0, 128}
\definecolor{pythonPurple}{RGB}{148,0,211}
\definecolor{darkgreen}{rgb}{0,0.3,0}
\definecolor{pythonOrange}{RGB}{255,165,0}

\definecolor{CBred}{RGB}{213, 94, 0}
\definecolor{CBblue}{RGB}{0,114,178}

\definecolor{VIColor}{HTML}{D55E00}
\definecolor{GVIColor}{HTML}{0072B2}
\definecolor{GVIColor1}{HTML}{0072B2}
\definecolor{GVIColor2}{HTML}{56B4E9}
\definecolor{FVIColor}{HTML}{009E73}
\definecolor{DVIColor}{HTML}{009E73}

\newcommand{\DVIColor}[1]{\textcolor{DVIColor}{#1}}%
\newcommand{\GVIColor}[1]{\textcolor{GVIColor1}{#1}}%
\newcommand{\VIColor}[1]{\textcolor{VIColor}{#1}}%
\newcommand{\lossColor}[1]{\textcolor{lossColor}{#1}}%
\newcommand{\PiColor}[1]{\textcolor{PiColor}{#1}}%
\newcommand{\DColor}[1]{\textcolor{DColor}{#1}}%
\newcommand{\Dcolor}[1]{\textcolor{DColor}{#1}}%

\usepackage{tikz} 
\usetikzlibrary{arrows,shapes}
\usepackage{bigints} 

\usepackage[utf8]{inputenc} 
\usepackage[T1]{fontenc}    
\usepackage{hyperref}       
\usepackage{url}            
\usepackage{booktabs}       
\usepackage{amsfonts}       
\usepackage{nicefrac}       
\usepackage{microtype}      
\usepackage{xspace}

\usepackage{xcolor} 
\usepackage[utf8]{inputenc}
\usepackage[english]{babel}
\usepackage{amsthm}
\usepackage{multicol}

\newcommand{\acro}[1]{\textsc{#1}\xspace}
\newcommand{\CP}{{\acro{\smaller CP}}}
\newcommand{\GP}{{\acro{\smaller GP}}}
\newcommand{\DGP}{\acro{\smaller DGP}}
\newcommand{\GBI}{{\acro{\smaller GBI}}}
\newcommand{\GVI}{\acro{\smaller GVI}}
\newcommand{\ELBO}{\acro{\smaller ELBO}}
\newcommand{\KLD}{\acro{\smaller KLD}}
\newcommand{\FVI}{\acro{\smaller F-VI}}
\newcommand{\Fsmall}{\acro{\smaller F}}

\newcommand{\BLR}{\acro{\smaller BLR}}
\newcommand{\BMM}{\acro{\smaller BMM}}

\newcommand{\renyi}{R\'enyi }

\newcommand{\TVD}{\acro{\smaller TVD}}
\newcommand{\AD}{\acro{\smaller AD$^{(\alpha)}$}}
\newcommand{\BD}{\acro{\smaller $D_{B}^{(\beta)}$}}
\newcommand{\GD}{\acro{\smaller $D_{G}^{(\gamma)}$}}
\newcommand{\ABGD}{\acro{\smaller $D_{G}^{(\alpha,\beta,r)}$}}
\newcommand{\RAD}{\acro{\smaller RD$^{(\alpha)}$}}
\newcommand{\wKLD}{\acro{\smaller$\frac{1}{w}$KLD}}
\newcommand{\RADfive}{\acro{\smaller $D_{AR}^{(0.5)}$}}

\newcommand{\ED}{\acro{\smaller ED}}
\newcommand{\JD}{\acro{\smaller JD}}
\newcommand{\FD}{\acro{\smaller FD}}

\newcommand{\ADqp}{\acro{{\smaller $D_{A}^{(\alpha)}$}$(q||\pi)$}}
\newcommand{\BDqp}{\acro{{\smaller $D_{B}^{(\beta)}$}$(q||\pi)$}}
\newcommand{\GDqp}{\acro{{\smaller $D_{G}^{(\gamma)}$}$(q||\pi)$}}
\newcommand{\RADqp}{\acro{{\smaller $D_{AR}^{(\alpha)}$}$(q||\pi)$}}

\newcommand{\paraBD}{\acro{\smaller$\*D_{B}^{(\beta)}$:}}
\newcommand{\paraGD}{\acro{\smaller$\*D_{G}^{(\gamma)}$:}}
\newcommand{\paraRAD}{\acro{\smaller$\*D_{AR}^{(\alpha)}$:}}
\newcommand{\parawKLD}{\acro{\smaller$\frac{1}{w}$KLD:}}

\newcommand{\VI}{\acro{\smaller VI}}
\newcommand{\ML}{\acro{\smaller ML}}
\newcommand{\EP}{\acro{\smaller EP}}
\newcommand{\BNN}{\acro{\smaller BNN}}
\newcommand{\RMSE}{\acro{\smaller RMSE}}
\newcommand{\MVN}{\acro{\smaller MVN}}
\newcommand{\SPD}{\acro{\smaller SPD}}
\newcommand{\VAEs}{\acro{\smaller VAEs}}
\newcommand{\VAE}{\acro{\smaller VAE}}
\newcommand{\bVAEs}{\acro{\smaller $\beta$-VAEs}}
\newcommand{\bVAE}{\acro{\smaller $\beta$-VAE}}

\newcommand{\VRBBVI}{\acro{\smaller VRBBVI}}
\newcommand{\BBVI}{\acro{\smaller BBVI}}
\newcommand{\VRBBGVI}{\acro{\smaller VRBBGVI}}
\newcommand{\BBGVI}{\acro{\smaller BBGVI}}

\newcommand{\LLN}{\acro{\smaller LLN}}

\newcommand{\Lb}{\acro{\smaller $\mathcal{L}_{p}^{\beta}$}}
\newcommand{\Lg}{\acro{\smaller $\mathcal{L}_{p}^{\gamma}$}}
\newcommand{\toL}{\longrightarrow}

\newcommand{\G}{$\Gamma$}

\newcommand{\jack}[1]{{{\color{red}Jack: [#1]}}}
\newcommand{\jeremias}[1]{{\color{blue}Jeremias: [#1]}}
\newcommand{\theo}[1]{{\color{orange}Theo: [#1]}}
\newcommand{\chris}[1]{{\color{purple}Chris: [#1]}}
\newcommand{\structure}[1]{{\color{darkgreen}structure: [#1]}}

\DeclareMathOperator*{\argmax}{arg\,max}
\DeclareMathOperator*{\argmin}{arg\,min}
\DeclareMathOperator*{\arginf}{arg\,inf}

\newcommand{\EPSRC}{\acro{\smaller EPSRC}}
\newcommand{\OxWaSP}{\acro{\smaller OxWaSP}}

\def\*#1{\bm{#1}} 
\newcommand{\f}[1]{\bm{#1}}%

\newenvironment{proofsketch}{%
  \renewcommand{\proofname}{Proof sketch}\proof}{\endproof}

\title{Frequentist Consistency of Generalized Variational Inference}

\author{
  Jeremias Knoblauch\\
  The Alan Turing Institute\\
  Dept. of Statistics\\
  University of Warwick\\
  \texttt{j.knoblauch@warwick.ac.uk}
}

\begin{document}

\maketitle

\begin{abstract} 
 This paper investigates Frequentist consistency properties of the posterior distributions constructed via Generalized Variational Inference (\GVI) \cite{GVI}.
 A number of generic and novel strategies are given for proving consistency, relying on the theory of $\Gamma$-convergence. 
 Specifically, this paper shows that under minimal regularity conditions, the sequence of \GVI posteriors is consistent and collapses to a point mass at the population-optimal parameter value as the number of observations goes to infinity.
 The results extend to the latent variable case without additional assumptions and hold under misspecification.
 Lastly, the paper explains how to apply the results to a selection of \GVI posteriors with especially popular variational families. For example, consistency is established for \GVI methods using the mean field normal variational family, normal mixtures, Gaussian process variational families as well as neural networks indexing a normal (mixture) distribution.
 %
\end{abstract}

\section{Introduction}

Bayesian inference methods are characterized by the desire to produce posterior belief distributions over a parameter $\*\theta \in \*\Theta$. 
The standard way of doing this is via Bayes' Rule: Given a prior belief $\pi(\*\theta)$ about the parameter and observations $x_{1:n}^o$ linked to $\*\theta$ via a likelihood function $p(x_i^o|\*\theta)$, the posterior is computed through a multiplicative updating rule as
\begin{IEEEeqnarray}{rCl}
    q^{\ast}(\*\theta) & \propto & \pi(\*\theta)\prod_{i=1}^n\exp\{-\ell(\*\theta, x_i^o)\},
    \label{eq:bayes_rule}
\end{IEEEeqnarray}
where $\ell(\*\theta, x_i^o) = -\log p(x_i^o|\*\theta)$.
While this way of writing Bayes rule might seem cumbersome, it reveals that the multiplicative structure is applicable to {any} loss function. In fact, replacing the negative log likelihood with any loss $\ell: \*\Theta \times \mathcal{X} \to \mathbb{R}$ for which the normalizer of eq. \eqref{eq:bayes_rule} exists provides a coherent and principled way to update beliefs about an arbitrary parameter $\*\theta$ \citep{Bissiri}.
For example, imagine that one wishes to update beliefs a robust measure of central tendency in the observation sequence $x_{1:n}^o$ in a Bayesian manner. 
A loss-based Bayesian treatment of this problem would combine a prior belief $\pi$ about the median $\*\theta^{\ast}$ with the loss $\ell(\*\theta, x_i^o) = |\*\theta - x_i^o|$. 
Together, these two ingredients yield the generalized Bayesian posterior given above, see for instance \citep{Jewson, RBOCPD, BissiriExp, RobustBayesGamma, GoshBasuPseudoPosterior} for some other interesting applications where $\ell(\*\theta, x_i^o) \neq -\log p(x_i^o|\*\theta)$.
%

Many theoretical properties make (generalized) Bayesian posteriors an attractive object to study for parameter inference.
One of these properties is consistency, which revolves around the posterior's asymptotic collapse to the population-optimum $\*\theta^{\ast}$ of $\*\theta$.
This property is well-studied for both standard \citep[e.g.][]{BayesianConsistencyClassic1, BayesConsistencyReview, BayesianConsistencyBVMMisspecification} and generalized \citep[e.g.][]{BayesianConsistencyGibbs,BayesianConsistencyMiller} posteriors. 
Contributions on this matter typically revolve around two interrelated questions: 
Firstly, one is interested whether the posterior concentrates at all \citep[see e.g.][]{BayesianConsistencyWalker, BayesianConsistencySafeBayes}. 
Secondly -- provided that the posterior does concentrate -- one is interested in the speed  at which the posterior contracts \citep[e.g.,][]{BayesianConsistencyAdaptive,BayesianConsistencyClassic2} or whether one can derive finite-sample PAC-bounds \citep[e.g.][]{BayesianConsistencyPACBayesClassic}.
While it is fair to say that Bayesian consistency holds under very mild regularity conditions so long as the parameter space is finite-dimensional, proving consistency can be far more demanding in the nonparametric case \citep[e.g.][]{BayesianConsistencyNonIID, BayesianNonparam1, BayesianNonparam2, BayesianNonparam3, BayesianConsistencyNonparam1, BayesianConsistencyNonparamMisspec, BayesianConsistencyJudith} and for Bayesian inverse problems \citep[e.g.][]{BayesianInverseBVM1, BayesianInverseBVM2, BayesianInverseBVM3}.

While beliefs derived from Bayes rule in eq. \eqref{eq:bayes_rule} can be written down analytically, any downstream computation (e.g., expectations with respect to $q^{\ast}(\*\theta)$) will require algorithms to sample from $q^{\ast}(\*\theta)$. 
Depending on the scale of the problem, this may become computationally infeasible. A popular approximation strategy geared towards alleviating the computational burden is Variational Inference (\VI) \citep[e.g.][]{JordanVI}. \VI approximates $q^{\ast}(\*\theta)$ with an element 
chosen from a parameterized set $\mathcal{Q}^{\*\theta}$ of distributions on $\*\Theta$.
%
%
At least three interpretations exist for justifying standard \VI: 
Firstly, one can see it as finding the $q \in \mathcal{Q}^{\*\theta}$ minimizing the Kullback-Leibler divergence \citep{kullback1951information} to $q^{\ast}(\*\theta)$. 
Secondly, one can interpret \VI as choosing the $q \in \mathcal{Q}^{\*\theta}$ maximizing a lower bound on the log evidence. 
Thirdly, relative to the variational characterization of the Bayesian posterior derived for instance by \citep{Csiszar}, \citet{DonskeerVaradhan} or \citet{Zellner}, \VI is the best possible $\mathcal{Q}^{\*\theta}$-constrained posterior belief about $\*\theta$ \citep{GVI}.

With increasing popularity of \VI, a substantial literature has evolved aiming at investigating whether the (approximate) posteriors produced by \VI are consistent. 
These efforts began early and usually focused on certain special cases of interest \citep[see e.g.][]{VBConsistencyExpModels, VBConsistencyMEstimator, VBConsistencyCommunityDetection} but have seen a substantial advancement in their applicability and generality over the past few years \citep[e.g.][]{VBConsistencyConvergenceRates}. 
The most notable consistency results for standard \VI posteriors establishing the speed of convergence and the consistency of mixtures are arguably the contributions of \citet{VBConsistencyYulong}, \citet{VBConsistencyWang} and \citet{VBConsistencyAlquier3}.
A number of novel results have also been established for \VI approximations of tempered posteriors \citep{VBConsistencyAlphaVI, VBConsistencyAlquier1} and generalized Bayesian posteriors \citep{VBConsistencyAlquier2}.

Lastly, a few select alternative approximation techniques have recently been shown to satisfy consistency properties under appropriately strict regularity conditions. For instance, variational approximations based on R\'enyi's $\alpha$-divergence can be shown to concentrate at the optimal parameter value \citep{VBConsistencyRenyiAlpha}. Note that this contains Expectation Propagaion as a special case, consistency properties of which were also studied under much more restrictive assumptions before \citep{VBConsistencyEP}.

The current paper adds to this last strand of literature by proving consistency of Generalized Variational Inference (\GVI), a novel method to produce posterior beliefs introduced in \citet{GVI}.
Since one of the most interesting features of \GVI is its modularity, the current paper derives proof strategies applicable to any modularly composed \GVI problem satisfying minimal assumptions.
While this yields very broadly applicable consistency results, it also means that these results are not strong enough to make statements about convergence speed. 
In this sense, the results are weaker than existing results for standard \VI procedures \citep[e.g.][]{VBConsistencyYulong, VBConsistencyWang}. They are however as strong as the few existing consistency results on alternative (non-standard) variational approximations, e.g. those based on R\'enyi's $\alpha$-divergence and Expectation Propagation in \citet{VBConsistencyRenyiAlpha}.

The remainder of the paper first provides a brief overview of both \GVI and a general road map used for proving consistency of \GVI posteriors. Second, notation and assumptions are explained. The main results about \GVI consistency follow in the third section. The emphasis is on tracing out the main arguments and sketching the big steps involved in the proofs, but full details on all results can be found in the Appendix. 
Lastly, the findings are illustrated on a number of simulation examples. These experiments confirm not only that consistency holds, but also that the speed of convergence will be a function of the prior regularizer $D$.
%
%





\section{Generalized Variational Inference (\GVI)}
\label{sec:GVI}


In what follows, the motivation behind \GVI and the Rule of Three decomposition is briefly recapitulated.
Next, a short high-level overview is given about the steps by which the current paper will succeed in proving consistency of \GVI posteriors. 

\subsection{Generalized Variational Inference (\GVI) and the Rule of Three}

Recently, \citet{GVI} evolved the logic of generalized Bayesian posteriors. Their generalized representation of Bayesian inference  takes inspiration from the variational representation of Bayesian inference  given by 
    \begin{IEEEeqnarray}{rCl}
    q^{\ast}(\*\theta) & = & \arginf_{q \in \mathcal{P}(\*\Theta)}\left\{
        \mathbb{E}_{q(\*\theta)}\left[
        \sum_{i=1}^n
            \ell(\*\theta, x_i^o)
        \right]
        +
        \KLD\left(q || \pi \right)
    \right\},  \;\; \KLD(q||\pi) = \mathbb{E}_{q(\*\theta)}\left[\log\dfrac{q(\*\theta)}{\pi(\*\theta)}\right], \;\; 
    \nonumber
\end{IEEEeqnarray}
which is for instance explored in \citet{Zellner}, but arguably much older \citep{Csiszar, DonskeerVaradhan}. 
%
%
%
Based on the variational formulation of Bayes' rule, \citep{GVI} provides an axiomatic derivation for a generalization of Bayesian inference beyond the loss-based posteriors in eq. \eqref{eq:bayes_rule}. 
%
Specifically, this generalization is an optimization problem over a subspace $\Pi^{\*\theta} \subset \mathcal{P}(\*\Theta)$, where $\mathcal{P}(\*\Theta)$ is the set of all probability measures $\mathcal{P}(\*\Theta)$ on $\*\Theta$. The optimization is then performed relative to a loss $\ell$ on the data and a divergence $D$ on the prior: 
\begin{IEEEeqnarray}{rCl}
    F_n(q) = F_n(q|x_{1:n}, \ell, D) & = & 
    \mathbb{E}_{q(\*\theta)}\left[
        \sum_{i=1}^n
            \ell(\*\theta, x_i^o)
        \right]
        +
        D\left(q || \pi \right)
        \label{eq:Fn_no_latent} \\
    q^{\ast}(\*\theta) & = & \arginf_{q \in \Pi^{\*\theta}}
        F_n(q). \nonumber
\end{IEEEeqnarray}

The current paper refines notation and logic introduced in \citep{GVI} by allowing for $\ell$ to additionally depend on latent variables $\*z_{i}$ on $\mathcal{Z}$.
Specifically, after $n$ observations have been made, the observations $x_{1:n}$ are allowed to depend on a finite and proportional number $S(n) = |\cup_{i=1}^ns(i)| \leq mn$ of latent variables $\*z_{1:S(n)}$. 
Here, $s(i)$ with $|s(i)| \leq m$ are the indices  of the latent variables entering $i$-th loss term together with $x_i$.
%
%
%
%
To notationally unify losses with and without latent variables, $\*x_i = (\*x_i^o, \*z_{s(i)})$ denotes the random variable generating $x_i = (x_i^o, z_{s(i)})$, so that $x_i$ is the $i$-th observation together with its relevant latent variables. 
%
Naturally, $\*x_i^o$ is the \textit{observed}   and $\*z_{s(i)}$ the \textit{latent} component.
%
%
%
%
%
%
Many inference problems fall into this broad modelling framework. The following list gives three examples that are both useful and popular in practice:
\begin{itemize}
    \item[(1)] 
    Observation-specific latent variable models, which encompasses Random Effects models. 
    This case corresponds to $s(i) = \{i\}$ and $\*x_i = (\*x_i^o, \*z_{i})$);
    \item[(2)]
    Time series, spatial and spatio-temporal settings where the indices of the latent space have a natural partial ordering.
    For instance, autoregressive time series settings with latent variables correspond to $s(i) = \{i-m, i-m+1, \dots, i\}$ and $\*x_i = (\*x_i^o, \*z_{(i-m):i})$.
    \item[(3)]
    A single set of global latent variables not depending on $n$, such as mixture models.
    This setting corresponds to $s(i) = \{1,2,\dots,m\}$, $S(n) = m$ for all $1\leq i \leq n$ and each $n \in \mathbb{N}$.
\end{itemize}
%
%
Suitably extending notation for latent variables  further, let $\Pi  \subseteq  \mathcal{P}(\*\Theta \times \mathcal{Z}^{S(n)})$ be a (sequence of) subspaces. 
Clearly, $\Pi$ depends on $n$, but this dependence is suppressed notationally for legibility.
Lastly, note that the theoretical analysis of the current paper focuses on the case where $\Pi$ is partitioned as $\Pi = \Pi^{\*\theta} \times \Pi_n^z$ for $\Pi^{\*\theta} \subseteq  \mathcal{P}(\*\Theta)$ and $\Pi^{z}_n \subseteq \mathcal{P}(\mathcal{Z}^{S(n)})$.
%

%
%
%

Finally putting everything together, one obtains a generalized Bayesian inference problem via the Rule of Three taking the form $P(\ell, D, \Pi)$.
Here, $P(\ell, D, \Pi)$ specifies a sequence of optimization problems indexed by $n$ whose three arguments are $\ell$, $D$ and $\Pi$ and which relate to one another via
\begin{IEEEeqnarray}{rCl}
    F_n(q,p)= F_n(q, p|x_{1:n}, \ell, D)  & = &
        \mathbb{E}_{q(\*\theta)}
        \left[
            \frac{1}{n}\sum_{i=1}^n\mathbb{E}_{p(\*z_{S(n)})}\left[\ell(\*\theta, x_i^o, \*z_{s(i)})\right]
        \right]  + \frac{1}{n}D\left(q||\pi\right),
    \label{eq:Fn_latent}
    \\
    q_n(\*\theta), q_n^z(\*z_{1:n}) & = & \arginf_{(q, p) \in (\Pi^{\*\theta}, \Pi^z_n)} F_n(q, p). 
    \nonumber 
\end{IEEEeqnarray}
Note that in the current paper the infimum replaces the minimum of the original generalization of the Bayesian inference problem in \citep{GVI}. This is done to ensure that (i) the problem is solvable for any $n \in \mathbb{N}$ and that (ii) the dirac delta at the population-optimum can be attained in the limit, even if it is not an element fo $\Pi$.
Moreover, notation in the remainder is simplified by suppressing the arguments $\ell$ and $D$ and abbreviating the objective as $F_n(q,p)$.

Whether one seeks to conduct inference with or without latent variables, the constituent parts of any generalized Bayesian inference problem are $\ell$, $D$ and $\Pi$. For an extensive discussion of their relationship, modularity and properties, see \citep{GVI}. 
Essentially, they determine the inference problem one wishes to solve. More specifically, they specify
\begin{itemize}
    \item A \textbf{loss} $\ell:\*\Theta \times \mathcal{X} \to \mathbb{R}$ linking a parameter of interest $\*\theta$ to the observations $x_{1:n}$.
    %
    Further, the loss is defined on $\mathcal{X} = \mathcal{X}^o \times \mathcal{Z}^m$, where $\mathcal{X}^o$ and $\mathcal{Z}^m$ denote the spaces of the observable random variables $\*x_i^o$ and the latent components $\*z_{s(i)}$, respectively.
    While the loss will be assumed to be additive over $x_{1:n}$, this is not the same as requiring independence between the random variables $\*x_{1:n}$ that generated the observations $x_{1:n}$.
    For instance, though $\ell(\*\theta, x_i)$ is allowed to depend on observations $x_j^o$ with $i\neq j$, this is suppressed from notation for easier legibility. 
    For example, for a loss relative to a time series with $p$-th order dependency, is viable to interpret $\ell(\*\theta, x_i)$ as $\ell(\*\theta, x_i; x_{i-p-1:i-1}^o)$.
    %
    %
    Notice also that while standard Bayesian inference requires the loss to be a negative log likelihood,  substantial advances have recently led to alternative loss functions amenable to principled Bayesian  inference \citep[e.g.][]{Bissiri,Jewson,AISTATSBetaDiv}.
    In line with this, the current paper does not limit itself to negative log likelihoods and derives results that hold for arbitrary losses.
    \item A {divergence ${D}:\mathcal{P}(\*\Theta)\times \mathcal{P}(\*\Theta) \to \mathbb{R}_+$ regularizing the parameter posterior} with respect to the prior $\pi$. As $D$ determines how
    uncertainty about the parameter $\*\theta$ is measured and quantified, it is called \textbf{uncertainty quantifier}. 
    To see why $D$ is called uncertainty quantifier, suppose that $D$ were absent from eq. \eqref{eq:Fn_no_latent},  which is to say that $D(q\|\pi) = 0$ for all $q \in \Pi^{\*\theta}$.
    Now, suppose further that $\Pi^{\*\theta}$ is big enough such that for all $\*\theta' \in \*\Theta$, one can find a sequence $q_k \in \Pi^{\*\theta}$ such that $q_k \overset{\mathcal{D}}{\rightarrow} \delta_{\*\theta'}$.
    In this case, it becomes clear that eq. \eqref{eq:Fn_no_latent} is solved by the point mass $\delta_{\widehat{\*\theta}_n}(\*\theta)$, where $\widehat{\*\theta}_n = \argmin_{\*\theta}\{ \sum_{i=1}^n \ell(\*\theta, x_i) \}$ is the empirical loss-minimizer.
    The effects of changing $D$ are manifold and can induce prior robustness as well as more conservative posterior marginals, see \citet{GVI}.
    %
    %
    \item A sequence of \textbf{admissible posteriors $\Pi \subseteq \mathcal{P}(\*\Theta)\times \mathcal{P}(\mathcal{Z}^{S(n)})$} which the regularized expected loss is minimized over.  
    Once again, dependence of this sequence on $n$ is  suppressed for legibility.
    In standard Bayesian inference, $ \mathcal{P}(\*\Theta\times\mathcal{Z}^{S(n)}) = \Pi$. For parameterized subsets of $ \mathcal{P}(\*\Theta)\times\mathcal{P}(\mathcal{Z}^{S(n)})$, denote $\mathcal{Q} = \mathcal{Q}^{\*\theta} \times \mathcal{Q}^z_n = \Pi$. Using this notation, $\mathcal{Q}$ could for instance denote the mean field normal variational family over $\*\Theta \times \mathcal{Z}^{S(n)}$.
\end{itemize} 
The Rule of Three can recover Bayes rule as in eq. \eqref{eq:bayes_rule} by taking the uncertainty quantifier to be the Kullback-Leibler divergence (\KLD) and the space of admissible posteriors is $\mathcal{P}(\*\Theta)$.
Using the abbreviation introduced above, this means that for a loss function $\ell$, Bayes rule solves $P(\ell, \KLD, \mathcal{P}(\*\Theta))$. 
Further, denoting $\mathcal{Q}^{\*\theta}$ as a a parameterized subset of $\mathcal{P}(\*\Theta)$, one may recover Variational Inference (\VI) from this form, too.
In other words, $P( \ell, \KLD, \mathcal{Q})$ corresponds to Variational Inference (\VI) relative to eq. \eqref{eq:bayes_rule}. 
%
%
Building on this insight, \citep{GVI} provides a principled motivation for tackling the more general class of problems with form $P(\ell, D, \mathcal{Q})$, calling them Generalized Variational Inference (\GVI) problems.
For the purpose of this paper, the original definition is extended to allow for latent variables, too.
{\begin{definition}[\GVI]
    {Any Bayesian inference method solving $P(\ell, D,  \mathcal{Q})$ with $\ell$ a loss function possibly depending on latent variables, $D$ a divergence and $\mathcal{Q}$ a parameterized subset of $\mathcal{P}(\*\Theta\times\mathcal{Z}^{S(n)})$ is a Generalized Variational Inference (\GVI) method.}
    \label{Def:GVI}
\end{definition}}

While \citep{GVI} clearly motivates advantages of \GVI over \VI, theoretical properties of such posteriors have remained unexplored until now.
The current paper provides the first step in filling this gap by proving frequentist consistency of \GVI procedures.
Since the most attractive feature of \GVI is its inherent modularity with respect to $\ell, D$ and $\mathcal{Q}$, the assumptions imposed to prove these consistency results are minimal.
Thus, the approach pursued here is a search for a minimal guarantee of \GVI, agnostic of the specific form of the problem apart from extremely weak regularity conditions on $\ell$, $D$ and  $\mathcal{Q}$.
%
%



%

\subsection{Proving \GVI consistency: A helicopter tour}

This paper sets out to provide a generic proof strategy applicable to \GVI with minimal assumptions. 
To this end, the proofs rely on functional and variational analysis. More specifically, the strategy deploys the machinery of $\Gamma$-convergence, which was for instance used in the work of \citet{VBConsistencyYulong} and \citet{VBConsistencyWang}.
\citet{gammaConvergence} provides a short, concise and rigorous introduction to the aspects of $\Gamma$-convergence relevant to the current work. 
%
%
Roughly speaking, the role of $\Gamma$-convergence in the present work is as follows: If a sequence of functions $F_n$ $\Gamma$-converges to a function $F$, then the sequence $q_n$ of its minimizers converges to the minimizer of $F$ under extremely mild regularity conditions.
To this end, the current paper studies the (stochastic) sequence of functions $F_n:\mathcal{Q}\to\mathbb{R}$  and the minimizers $q_n$ as given by eqs. \eqref{eq:Fn_latent} and \eqref{eq:Fn_no_latent}.
%
Because intuitions and proofs are the same with or without latent variables, it is instructive to keep things as simple as possible. 
To this end, the following section uses the case without latent variables in eq. \eqref{eq:Fn_no_latent} to explain this paper's proof strategies. 
Naturally, only parameterized subsets $\Pi^{\*\theta} = \mathcal{Q}^{\*\theta}$ of $\mathcal{P}(\*\Theta)$ are considered from here on out.
Before proceeding, one may wish to take another careful look at the optimization problem. 
%
\begin{IEEEeqnarray}{rCl}
    F_n(q) 
    & = & 
    \mathbb{E}_{q}\left[
        \frac{1}{n}\sum_{i=1}^n\ell(\*\theta, x_i)
    \right]
    + 
    \frac{1}{n}D(q||\pi)
    \nonumber \\
    q_n & = & \arginf_{q \in \mathcal{Q}^{\*\theta}}F_n(q).
    \nonumber
\end{IEEEeqnarray}
First, note that  mild regularity conditions ensure that $\frac{1}{n}\sum_{i=1}^n\ell(\*\theta, x_i) \overset{\mu-a.s.}{\toL} \mathbb{E}_{\mu}\left[
            \ell(\*\theta, \*x)
\right]$ as $n\to\infty$ for an appropriate 
probability measure $\mu$ on $\mathcal{X}$. 
Accordingly, it usually holds that $\widehat{\*\theta}_n = \argmin\{\frac{1}{n}\sum_{i=1}^n\ell(\*\theta, x_i)\} \to \*\theta^{\ast} = \argmin\{\mathbb{E}_{\mu}\left[
            \ell(\*\theta, \*x)
\right]\}$ as $n\to \infty$, $\mu$-almost surely.
%
Intuitively then, one expects the sequence $q_n$ to converge in distribution to $\delta_{\*\theta^{\ast}}(\*\theta)$ under mild regularity conditions, $\mu$-almost surely. 
In other words, the remainder of the paper will show that for $\mathcal{F}$ the set of continuous, bounded functions from $\*\Theta$ to $\mathbb{R}$,
\begin{IEEEeqnarray}{rCl}
    \mathbb{P}_{\mu}\left( q_n(\*\theta) \overset{\mathcal{D}}{\to} \delta_{\*\theta^{\ast}}(\*\theta) \right)  = 
     \mathbb{P}_{\mu}\left(
        \forall f \in \mathcal{F}:
        \lim_{n\to\infty}\left\{\mathbb{E}_{q_n}\left[f(\*\theta)\right] - f(\*\theta^{\ast})\right\} = 0
     \right)
    & = & 1.
    \label{eq:weak_convergence_almost_surely}
\end{IEEEeqnarray}
For certain cases, this paper also derives the weaker result where the above holds in the $\mu$-probability limit. Saying that the sequence $q_n$ weakly converges to $\delta_{\*\theta^{\ast}}(\*\theta)$ in the $\mu$-probability limit means that 
\begin{IEEEeqnarray}{rCl}
    \forall f \in \mathcal{F}, \lim_{n\to\infty}\mathbb{P}_{\mu}\left(
        \mathbb{E}_{q_n}\left[f(\*\theta)\right] - f(\*\theta^{\ast}) = 0
     \right)
    & = & 1.
    \label{eq:weak_convergence_in_P}
\end{IEEEeqnarray}
However, any direct way of showing that this intuition holds would require establishing $\Gamma$-convergence of $F_n$ to $\mathbb{E}_{q}\left[
        \mathbb{E}_{\mu}\left[
            \ell(\*\theta, \*x)
        \right]
    \right]$.
Unfortunately, the stochasticity of $F_n$ makes it hard to prove this result directly. 
The key insight of the current paper is to circumvent these technical complications associated with $F_n$'s stochasticity by instead analyzing the \textit{deterministic} sequence of functions $\overline{F}_n$ and its minimizers $\overline{q}_n$ given by 
\begin{IEEEeqnarray}{rCl}
    \overline{F}_n(q)
    & = &
    \mathbb{E}_{q}\left[
        \mathbb{E}_{\mu}\left[
            \ell(\*\theta, \*x)
        \right]
    \right]
    + 
    \frac{1}{n}D(q||\pi)
    \nonumber \\
    \overline{q}_n & = & \arginf_{q \in \mathcal{Q}^{\*\theta}}\overline{F}_n(q).
    \nonumber
\end{IEEEeqnarray}
For this new auxiliary objective $\overline{F}_n$, establishing the desired  $\Gamma$-convergence too
to $\mathbb{E}_{q}\left[
        \mathbb{E}_{\mu}\left[
            \ell(\*\theta, \*x)
        \right]
    \right]$ 
turns out to be much simpler.
The last and most important part of the proofs will then be to show that the sequences $q_{n}$ and $\overline{q}_n$ become arbitrarily close as $n\to\infty$ ($\mu$-almost surely or in the $\mu$-probability limit).
More precisely, this paper shows that the sequence $\{q_{n}\}_{n=1}^{\infty}$ constitutes a sequence of $\varepsilon_n$-minimizers of $\overline{F}_n$, i.e.
\begin{IEEEeqnarray}{rCl}
    \overline{F}_n(q_n) \leq \inf_{q \in \mathcal{Q}^{\*\theta}} \overline{F}_n(q) + \varepsilon_n,
    \nonumber
\end{IEEEeqnarray}
where $\varepsilon_n$ is a stochastic sequence converging to zero ($\mu$-almost surely or in the $\mu$-probability limit).
This -- together with $\Gamma$-convergence and equi-coerciveness of $\overline{F}_n$ -- suffices to show that as desired, eq. \eqref{eq:weak_convergence_almost_surely} holds.
%
%
To summarize, this paper shows consistency of \GVI in three steps:
\begin{itemize}
    \item[(1)]
    Establishing that
    $\overline{F}_n$ is equi-coercive and $\Gamma$-converges to $\mathbb{E}_{q}\left[
        \mathbb{E}_{\mu}\left[
            \ell(\*\theta, \*x)
        \right]
    \right]$, from which it follows that $\overline{q}_n \overset{\mathcal{D}}{\to} \delta_{\*\theta^{\ast}}$ as $n\to\infty$, i.e. the minimizers of $\overline{F}_n$ converge;
    \item[(2)]
    Showing that the minimizers $q_n$ of the stochastic sequence ${F}_n$ are $\varepsilon_n$-minimizers of $\overline{F}_n$;
    \item[(3)]
    Proving that $\varepsilon_n$ goes to zero $\mu$-almost surely (in the $\mu$-probability limit) as $n\to\infty$. This together with the first two findings finally implies that as desired, $q_n \overset{\mathcal{D}}{\to} \delta_{\*\theta^{\ast}}$ holds $\mu$-almost surely (in the $\mu$-probability limit).
\end{itemize}
Perhaps surprisingly, these are the exact same steps one follows for the latent variable case, too. 
In this case, one essentially treats the additional posterior $q_n^z$ over the latent variables in the same way one treats nuisance parameters.



\section{Preliminaries}



Throughout, $(\*\Theta, \|\cdot\|)$ will be  a normed space of finite dimension. Notice that normed spaces are always metrizeable. Most of the derived results hold more generally: In fact, the metric space property is only important for Lemma \ref{lemma:equi_coercive_Fbar} and could be weakened if so desired. 
%
%
%

To make the theoretical developments concise, certain notational liberties are taken: Suppose the measure $\nu \in \mathcal{P}(\*\Theta)$ admits a density $q_{\nu}$ on $\*\Theta$ that is absolutely continuous with respect to the Lebesgue measure. 
The current paper treats the density  $q_{\nu}$ as if it were an element of $\mathcal{P}(\*\Theta)$ itself. This circumvents the need to carefully distinguish measures (which are elements of $\mathcal{P}(\*\Theta)$) and their densities (which are not). 
Though this greatly simplifies conceptualizing the key problems and proofs of the paper, it also leads to slightly abusive notation: For example, whenever one writes $q_n \overset{\mathcal{D}}{\to} \delta_{\*\theta^{\ast}}$, the statement one really makes is that the sequence of measures $\nu_{n} \in \mathcal{P}(\*\Theta)$ with densities $q_n$ converges weakly to the measure $\delta_{\*\theta^{\ast}} \in \mathcal{P}(\*\Theta)$. %


Further, it is assumed that the fixed numbers $x_{i} \in \mathcal{X}$ are realizations of random variables $\*x_i:\Omega \to \mathcal{X}$ defined on potentially different probability spaces $({\Omega}, \mathcal{F}, \mu_i)$.
%
This implies that $\*x_{1:n}$ is allowed to exhibit dependency. In other words, the derived results do not require $\*x_{1:n}$ to be independent and identically distributed (iid) copies of the random variable $\*x_1$.
Instead and as Assumption \ref{AS:min_exists} will show, only a strong law of large numbers needs to hold. Specifically, one needs that for some probability measure $\mu$ on the measureable space $(\Omega, \mathcal{F})$, $\frac{1}{n}\sum_{i=1}^n\ell(\*\theta, x_i) \overset{\mu-a.s.}{\toL} \mathbb{E}_{\mu}[\ell(\*\theta, \*x)]$ as $n\to \infty$.
While $\mu$ has a straightforward interpretation as the probability measure of $\*x_1$ when $x_i \overset{iid}{\sim} \*x_1$, any probability measure on $(\Omega, \mathcal{F})$ is sufficient, see also Remark \ref{remark:AS1} for an example.

Lastly, note that in the presence of latent variables,  the loss is defined on $\mathcal{X} = \mathcal{X}^o \times \mathcal{Z}^m$, where $\mathcal{X}^o$ and $\mathcal{Z}^m$ denote the spaces of the observables $x_i^o$ and the latent components $z_{s(i)}$.
Similarly, in the absence of latent variables, it holds that $\mathcal{X} = \mathcal{X}^o$.


\section{Assumptions}
\label{sec:assumptions}

%



In the following, the assumptions used for proving \GVI consistency are explained. 
In a nutshell, Assumptions \ref{AS:min_exists}, \ref{AS:dirac}, \ref{AS:D}, \ref{AS:suitable} and \ref{AS:finite_solution_exists} are mild regularity conditions that hold in practice for virtually any \GVI problem of interest.
Assumption \ref{AS:varepsilon_convergence} is less harmless and will have to be verified. However and as Section \ref{sec:varepsilon-convergence} investigates in detail, there is a plethora of interpretable conditions that imply Assumption \ref{AS:varepsilon_convergence}. 
For example, Assumption \ref{AS:varepsilon_convergence}  holds if (i) the prior $\pi(\*\theta)$ is normal and $\mathcal{Q}$ is a normal mean-field variational family (or any family encompassing it), 
and (ii) all $x_i$ are independent copies of the same random variable $\*x_1$.

%
%


\begin{assumption}
    The \GVI problem $P( \ell, D,\mathcal{Q})$ is well-defined:  It holds that
    \begin{itemize}
        %
         \item[(1)]
            The loss function $\ell:\*\Theta \times \mathcal{X} \to \mathbb{R}$ is discontinuous at most at finitely many points.
        \item[(2)] For any $x_i^o \in \mathcal{X}^o$ and any $n$,  $\ell(\*\theta, x_i) = \ell(\*\theta, x_i^o, z_{s(i)}) < \infty$ for all $z_{s(i)} \in \mathcal{Z}^m$.
         \item[(3)]
         The minimizers $\hat{\*\theta}_n  = \argmin_{\*\theta}\left\{\frac{1}{n}\sum_{i=1}^n\ell(\*\theta, x_i)\right\} \in \*\Theta$ exist for all $n$;
    \end{itemize}
        Moreover, the data-generating mechanism is well-behaved: For a probability measure $\mu$ on $\mathcal{X}$,
    \begin{itemize}
         \item[(4)]
        The loss satisfies a law of large numbers, i.e.
        $\frac{1}{n}\sum_{i=1}^n\ell(\*\theta, x_i) \overset{\mu-a.s.}{\toL} \mathbb{E}_{\mu}\left[ \ell(\*\theta, \*x) \right]$;
         \item[(5)]
        The $\mu$-population-minimizer ${\*\theta}^{\ast} = \argmin_{\*\theta}\mathbb{E}_{\mu}\left[ \ell(\*\theta, \*x) \right] \in \*\Theta$  exists and is unique;
        \item[(6)]
          The loss is finite in $\mu$-expectation, i.e. $\mathbb{E}_{\mu}\left[ \ell(\*\theta, \*x) \right]<\infty$ for all $\*\theta \in \*\Theta$;
        
        \item[(7)]  One of the following holds true: 
            $\mathbb{E}_{\mu}\left[ \ell(\*\theta, \*x) \right]$ is coercive in $\*\theta$ or
            $\*\Theta$ is compact.
    \end{itemize}
    \label{AS:min_exists}
\end{assumption}
\begin{remark}
    This set of assumptions is extremely mild and can be assumed to hold in virtually all situations of interest.
    For example, the requirement in (7) that $\mathbb{E}_{\mu}\left[ \ell(\*\theta, \*x) \right]$ be coercive is strictly more general than assuming convexity in $\*\theta$.
    Naturally, this comes at a cost: Unlike the stronger assumptions typically imposed for proving consistency, the above is not enough to establish concentration rates. 
    %
    %
    %
    \label{remark:AS1_mild}
\end{remark}

\begin{remark}
    The interpretation of $\mu$ for the case where $x_i \overset{iid}{\sim} \*x_1$ is clearly that of the probability measure corresponding to $\*x_1$.
    Things are perhaps less obvious in the dependent case, albeit not conceptually complicated. 
    For example, suppose that (i) $\ell(\*\theta, x_i) = \ell(\*\theta, x_i^o) = \ell(\*\theta, x_i^o; x_{i-1}^o)$ is the likelihood of a first order autoregressive process without latent variables and (ii) this autoregression accurately describes how $x_{1:n}$ was generated.  Then -- provided that a strong law of large numbers holds -- $\mu$ is a conditional probability measure on the events $\{\*x_i = x_i|\*x_{i-1} = x_{i-1}\}$.
    Dependencies like these are notationally suppressed for readability, but do not affect any of the results derived in the current paper.
    %
    %
    %
    %
    \label{remark:AS1}
\end{remark}
\begin{remark}
    In addition to dependence, one may well be interested in the convergence properties of \GVI posteriors built with a sequence of heterogeneous losses $\{\ell_i(\*\theta, x_i)\}_{i=1}^{\infty}$ where $\ell_i \neq \ell_j$ for some $i,j$.
    In this case, all derived convergence results follow after an easy adaption of the the above assumption.
    Specifically, one requires that (1) and (2) hold for each loss $\ell_i$ and that the minimizers in (3) exist for $\frac{1}{n}\sum_{i=1}^n\ell_i(\*\theta, x_i)$ instead. 
    Further, one requires that there exists some function $\widetilde{\ell}:\*\Theta \times \mathcal{X} \to \mathbb{R}$ such that $\frac{1}{n}\sum_{i=1}^n\ell_i(\*\theta, x_i) \overset{\mu-a.s.}{\longrightarrow} \mathbb{E}_{\mu}\left[\widetilde{\ell}(\*\theta, \*x)\right]$.
    Replacing the old convergence requirement in (4) with the new one and $\mathbb{E}_{\mu}\left[\ell(\*\theta, \*x)\right]$ with $\mathbb{E}_{\mu}\left[\widetilde{\ell}(\*\theta, \*x)\right]$ in (5), (6) and (7) completes the adaption to heterogeneous losses.
    This adapted assumption can now be directly used without any other additional requirements to replace the original Assumption \ref{AS:min_exists} for \textit{all} theoretical results derived in the sequel.
    %
    %
\end{remark}

The following assumption makes sure that the space $\mathcal{Q}$ over which one solves the \GVI problem admits convergence to a point mass. This is satisfied by virtually all approximate posteriors used in practice. Perhaps most importantly, mean field normal distributions satisfy this requirement.
%
%
\begin{assumption}
    The variational family $\mathcal{Q} = \mathcal{Q}^{\*\theta} \times \mathcal{Q}_n^z$ with $ \mathcal{Q}^{\*\theta} = \{q(\*\theta|\*\kappa): \*\kappa \in \*K \}$ and $\mathcal{Q}_n^z = \{q^z_n(\*z_{S(n)}|\*\eta): \*\eta \in \*H \}$ consists of absolutely continuous densities with respect to the Lebesgue measure. Moreover, for all $n \in \mathbb{N}$ and any $(\*\theta^{\ast},z) \in \*\Theta \times \mathcal{Z}^{S(n)}$, there exist sequences $\{\*\kappa_k\}_{k=1}^{\infty}$ and $\{\*\eta_k\}_{k=1}^{\infty}$ of variational parameters so that $ q(\*\theta|\*\kappa_k) \overset{\mathcal{D}}{\to} \delta_{{\*\theta}^{\ast}}(\*\theta)$ and $ q^z_n(\*z_{S(n)}|\*\eta_k) \overset{\mathcal{D}}{\to} \delta_{z}(\*z_{S(n)})$ as $k\to\infty$.
    \label{AS:dirac}
\end{assumption}
%

Beyond the variational family, one needs to ensure some minimal regularity properties of $D$. These properties are satisfied by virtually any statistical divergence, including all $f$- and Bregman divergences as well as the family of $\alpha\beta\gamma$-divergences \citep{ABCdiv}. 
%
\begin{assumption}
    The \GVI uncertainty quantifier $D:\mathcal{P}(\*\Theta)^2 \to \mathbb{R}_{+}$ is a statistical divergence. Further, it is lower semi-continuous in its first argument with respect to the weak topology of $\mathcal{P}(\*\Theta)$.
    \label{AS:D}
\end{assumption}

Next, one needs to ascertain that prior and uncertainty quantifier $D$ are suitable for the variational family $\mathcal{Q}$. While this may look like a strong requirement at first glance, it is satisfied in practice for all but the most pathological situations. 
%
%
\begin{assumption}
    The prior $\pi$ and the \GVI uncertainty quantifier $D$ are suitable for the variational family $\mathcal{Q}^{\*\theta}$: For all $q\in\mathcal{Q}^{\*\theta}$, $D(q||\pi) < \infty$.
    \label{AS:suitable}
\end{assumption}
\begin{remark}
    Note that the function of this assumption is similar to the requirement that $\pi(\*\theta) > 0$ in a neighbourhood of $\*\theta^{\ast}$ in traditional Bayesian consistency proofs.
    Its role is to ensure that concentration on the dirac measure at $\*\theta^{\ast}$ is possible. In fact, the requirement is slightly stronger than necessary for consistency: It would suffice to require that there exists a sequence $p_n \in \mathcal{Q}^{\*\theta}$ so that (i) $p_n \overset{\mathcal{D}}{\to} \delta_{\*\theta^{\ast}}$ and (ii) $D(p_n||\pi) < \infty$ for all finite $n$. 
    If $D = \KLD$, it is clear that this latter requirement is satisfied if and only if $\mathcal{Q}$ satisfies Assumption \ref{AS:dirac} and $\pi(\*\theta) > 0$ in a neighbourhood of $\*\theta^{\ast}$.
    %
    %
    Indeed, this equivalence holds over the wider class of choices for $D$ in the set 
    \begin{IEEEeqnarray}{rCl}
        \{D: \text{$D$ satisfies Assumption \ref{AS:D} and }D(q||\pi) = \infty \text{ $\forall q$ not absolutely continuous w.r.t. $\pi$} \}.
        \nonumber
    \end{IEEEeqnarray} 
    %
    Examples of uncertainty quantifiers in this set are the \KLD, the $\alpha$-divergence, R\'enyi's $\alpha$-divergence as well as the family of $f$-divergences and most Bregman divergences of practical interest.
\end{remark}

Lastly, one needs to ensure that the posteriors (i) do not represent infinitely bad beliefs and (ii) are not worse than the prior beliefs. 
%
%
%
The following assumption is an easy and intuitively appealing way to operationalize this requirement. 
%
\begin{assumption}
    The prior belief $\pi$ about $\*\theta$ is not infinitely bad: $\mathbb{E}_{\pi}\left[ \mathbb{E}_{\mu}\left[ \ell(\*\theta, \*x) \right] \right] = C_{\pi} < \infty$ 
    %
    Moreover, $\mathcal{Q}^{\*\theta}$ contains the singleton $\pi(\*\theta)$. In other words, $\mathcal{Q}^{\*\theta} = \{q(\*\theta|\*\kappa): \*\kappa \in \*K \}\cup\{\pi(\*\theta) \}$.
    %
    %
    \label{AS:finite_solution_exists}
\end{assumption}
\begin{remark}
    %
    %
    %
    The first part of this assumption has an interpretation as quality requirement on the prior: Having the belief $\pi$ about $\*\theta$ should not yield an infinite expected loss.
    While hard to verify, it is natural to assume that that this prior quality assumption holds in all situations of practical interest.
    The second part of the assumption is purely technical and makes the proofs less complicated, but has no bearings on inferential practice. 
    For the proofs, its role is simply to ensure that the \GVI posteriors do not become \textit{worse} than the prior belief after observing data. 
    Since this is not a situation that ever occurs in practice, one may safely remove $\pi$ from $\mathcal{Q}^{\*\theta}$ when using \GVI on a real-world problem.
    %
    %
    %
    %
    Functionally, the role of Assumption \ref{AS:finite_solution_exists} is to guarantee that there exist $q_n', \overline{q}_n' \in \mathcal{Q}^{\*\theta}$ such that $F_n(q_n') < \infty$ ($\mu$-a.s.) and $\overline{F}_n(\overline{q}_n')<\infty$ for all $n$.
    Similarly, note that if Assumptions \ref{AS:dirac} and \ref{AS:min_exists} also hold, then $q_n^z(\*z_{S(n)}) = \delta_{z_{S(n)}}(\*z_{S(n)})$ is an admissible choice for which $F_n(q_n', \delta_{z_{S(n)}}) <\infty$, too. 
    %
    %
    %
    In other words, the remarkably harmless assumption above guarantees that the sequence of \GVI problems produces posteriors that are not infinitely bad relative to both the population as well as the actually observed data.
    %
    %
\end{remark}


The next assumption is used to ensure that the parameter posterior vanishes everywhere on $\*\Theta$ except around the optimum value $\*\theta^{\ast}$. 
The role of this assumption as well as various sufficient conditions ensuring that it holds are discussed in section \ref{sec:varepsilon-convergence}.

\begin{assumption}
    For
    $\varepsilon_n = \int_{\*\Theta}\left[\frac{1}{n}\sum_{i=1}^n\ell(\*\theta, x_i) - \mathbb{E}_{\mu}\left[\ell(\*\theta, \*x)\right]\right] \overline{q}_n(\*\theta)d\*\theta$,
    one of the following holds:
    \begin{itemize}
        \item[(A)] $\varepsilon_n \overset{\mu-\mathcal{P}}{\toL} 0$ as $n\to \infty$;
        \item[(B)] $\varepsilon_n \overset{\mu-a.s.}{\toL} 0$ as $n\to \infty$.
    \end{itemize}
    \label{AS:varepsilon_convergence}
\end{assumption}
\begin{remark}
    Clearly, this is \textit{not} a harmless assumption per se. 
    Yet as this paper shows in section \ref{sec:varepsilon-convergence}, various constellations of mild assumptions suffice to prove that $\varepsilon_n$ converges to zero for rather different cases.
    For example, Assumption \ref{AS:varepsilon_convergence}(B) holds if $x_{1:n} \overset{iid}{\sim} \*x_1$, $\pi$ is a normal distribution and $\mathcal{Q}$ a mean-field normal family. 
    Further, it trivially holds on compact parameter spaces $\*\Theta$ or if $\ell(\*\theta, x_i) \leq M$ for some $M$ (as is often the case with losses of robust estimators).
    %
    %
\end{remark}

While this section does not discuss Assumption \ref{AS:varepsilon_convergence} further, the following is an auxiliary supposition useful in proving that Assumption \ref{AS:varepsilon_convergence} holds for a large class of \GVI problems.

\begin{assumption}
    There exists a compact subset $A\subset \*\Theta$ so that (i) $\*\theta^{\ast} \in A$ and
    (ii) $\pi\geq \overline{q}_n$  on $\*\Theta\setminus A$, for all $n \geq N$ for some $N<\infty$. 
    \label{AS:A_exists}
\end{assumption}
\begin{remark}
    As Corollary  \ref{corollary:consistency_Fbar} will show, $\overline{q}_n \overset{D}{\to} \delta_{\*\theta^{\ast}}$.
    Consequently, in most cases $\pi$ will eventually dominate $\overline{q}_n$ outside of a set $A$ containing $\*\theta^{\ast}$.
    Conceptually, this is a very mild requirement: It essentially says that $\pi$ has to be \textit{structurally} at most as concentrated as the elements in $\mathcal{Q}^{\*\theta}$. As both $\pi$ and $q \in \mathcal{Q}^{\*\theta}$ are parameterized, this structure is condensed into their parameter spaces\footnote{
        In fact, it is appropriate to think about this notion of structural concentration in terms of the statistical manifolds of $\pi$ and  $\mathcal{Q}^{\*\theta}$ (which are invariant to re-parameterizations) rather than their particular parameterizations.
    }.
    For instance, if $\pi$ is a normal distribution and $\mathcal{Q}^{\*\theta}$ is the mean field normal variational family, such a set $A$ can be found (Lemma \ref{lemma:A_exists_for_MFN}).
    %
    \label{remark:A_exists}
\end{remark}


\section{Main Results}

This section states the main results and sketches important aspects of the proofs. Details and a rigorous treatment of the corresponding results are discussed in Appendix \ref{sec:proof_details}.
In fact, this paper's findings permit substantially stronger results than the ones presented in the current section at the expense of verifying additional conditions.
In this light, the following Theorems are summaries of the most poignant and practically relevant findings derived in Appendix \ref{sec:proof_details}.
%
%
%

\begin{theorem}[Generic \GVI consistency]
    If Assumptions \ref{AS:min_exists}, \ref{AS:dirac}, \ref{AS:D}, \ref{AS:suitable}, \ref{AS:finite_solution_exists} and \ref{AS:varepsilon_convergence} hold, then \GVI posteriors are consistent. If Assumption \ref{AS:varepsilon_convergence}(A) holds, they are consistent in the $\mu$-probability limit. If Assumption \ref{AS:varepsilon_convergence}(B) holds instead, they are consistent  $\mu$-almost surely. 
    \label{thm:mainres_sketch}
\end{theorem}
\begin{proofsketch}
The three steps of proving \GVI consistency are outlined in three separate parts:
\begin{itemize}
    \item[]\textbf{Appendix \ref{sec:gamma_convergence}}:
    By virtue of Assumptions \ref{AS:min_exists}, \ref{AS:dirac}, \ref{AS:D}, \ref{AS:suitable}, two things can be shown to hold true. Firstly, one can show that the sequence of functions $\overline{F}_n$ is equi-coercive (Lemmas \ref{lemma:coercive_Psi} and \ref{lemma:equi_coercive_Fbar}). Secondly, one can show that it $\Gamma$-converges to $\mathbb{E}_{q}\left[\mathbb{E}_{\mu}\left[\ell(\*\theta,\*x) \right]\right]$ (Lemma \ref{lemma:gamma_conv_Fbar}). 
    Together, this implies that one of the main workhorses for proving consistency of \GVI can be deployed. Specifically, Corollary 7.24 in \citep{gammaConvergence} holds, proving that $\overline{q}_n \overset{\mathcal{D}}{\to} \delta_{\*\theta^{\ast}}$ (Corollay \ref{corollary:consistency_Fbar}). 
    \item[]\textbf{Appendix \ref{sec:varepsilon-minimizers}}:
    As Assumption \ref{AS:finite_solution_exists} holds, one can ensure that (i) the prior belief is not infinitely bad and (ii) posterior beliefs improve upon the prior. 
    This is formalized in Lemma \ref{lemma:Fubini}, together with the insight that $\ell(\*\theta, x_i) < \infty$ holds $\mu$-almost surely.
    Together with Assumptions \ref{AS:min_exists}, \ref{AS:dirac}, \ref{AS:D}, \ref{AS:suitable}, this suffices to prove that $\{q_n\}_{n=1}^{\infty}$ is a sequence of $\varepsilon_n$-solutions of $\overline{F}_n$ (Lemma \ref{lemma:eps_minimizers}). In other words, it holds that
    \begin{IEEEeqnarray}{rCl}
        \overline{F}_n(q_n) 
        & \leq &
        \inf_{q\in\mathcal{Q}^{\*\theta}}\overline{F}_n(q) + \varepsilon_n
        \nonumber
    \end{IEEEeqnarray}
    for a $\mu$-almost surely finite-valued sequence $\{\varepsilon_n\}_{n=1}^{\infty}$ whose terms are given by
    \begin{IEEEeqnarray}{rCl}
        \varepsilon_n & = &
        2\left| 
           \mathbb{E}_{\overline{q}_n}\left[
                \frac{1}{n}\sum_{i=1}^n\ell(\*\theta, x_i) - \mathbb{E}_{\mu}\left[\ell(\*\theta, \*x)\right]
           \right]
        \right|.
        \nonumber
    \end{IEEEeqnarray}

    \item[]\textbf{Appendix \ref{sec:varepsilon-convergence}}: 
    Having established the form of $\varepsilon_n$, the next step consists in deriving sufficient conditions for Assumption \ref{AS:varepsilon_convergence} to hold.
    To this end, recall that 
    \textbf{Appendix \ref{sec:gamma_convergence}} showed that $\overline{q}_n \overset{\mathcal{D}}{\to} \delta_{\*\theta^{\ast}}$ $\mu$-almost surely  (in the $\mu$-probability limit).
    From this, it immediately becomes clear that $\mathbb{E}_{\overline{q}_n}\left[\mathbb{E}_{\mu}\left[\ell(\*\theta^{\ast}, \*x)\right]\right] \longrightarrow \mathbb{E}_{\mu}\left[\ell(\*\theta^{\ast}, \*x)\right]$ $\mu$-almost surely  (in the $\mu$-probability limit).
    It follows that $\varepsilon_n \toL 0$ $\mu$-almost surely  (in the $\mu$-probability limit) as $n\to \infty$ only if 
    \begin{IEEEeqnarray}{rCl}
        \frac{1}{n}\sum_{i=1}^n\mathbb{E}_{\overline{q}_n}\left[
                \ell(\*\theta, x_i)\right]
                & \longrightarrow &
                \mathbb{E}_{\mu}\left[\ell(\*\theta^{\ast}, \*x)\right]
                \nonumber
    \end{IEEEeqnarray}
    $\mu$-almost surely (in the $\mu$-probability limit).
    Notice also that showing the above is equivalent to proving a law of large numbers for the triangular array 
    \begin{IEEEeqnarray}{rCl}
        \left\{ \left\{ \mathbb{E}_{\overline{q}_n}\left[ 
        \ell(\*\theta, x_i)
        \right] \right\}_{i=1}^n \right\}_{n=1}^{\infty}
        \nonumber
    \end{IEEEeqnarray}
    As it turns out, this can be shown to hold true under a variety of conditions that fall into two broad categories: Sufficient conditions for a dominated convergence theorem (e.g., Corollary \ref{corollary:summary_consistency}) as well as boundedness conditions on $\ell$ (e.g., Corollary \ref{corollary:summary_consistency_dependence} and Theorem \ref{thm:consistency_dependence_sketch}). 
    Provided that the loss is unbounded, the most practical way to establish conditions for dominated convergence is by verifying Assumption \ref{AS:A_exists} (e.g., Lemma \ref{lemma:A_exists_for_MFN}).
    %
    Finally, invoking the main workhorse --  Corollary 7.24 in \citep{gammaConvergence} -- once again, this implies that as desired, $q_n\overset{\mathcal{D}}{\to} \delta_{\*\theta^{\ast}}$ holds $\mu$-almost surely (in the $\mu$-probability limit).
\end{itemize}
As desired, this finally implies consistency of \GVI parameter posteriors.
\end{proofsketch}

As discussed Section \ref{sec:assumptions}, Assumptions \ref{AS:min_exists}, \ref{AS:dirac}, \ref{AS:D}, \ref{AS:suitable}, \ref{AS:finite_solution_exists} are all benign. Thus, the results in Appendix \ref{sec:gamma_convergence} and  \ref{sec:varepsilon-minimizers} outlined in the proof sketch above will almost always apply. 
The same cannot be said about the third step in Appendix \ref{sec:varepsilon-convergence}, which derives conditions under which Assumption \ref{AS:varepsilon_convergence} and thus consistency holds.
%
%
%
The next result gives an attractive way of instead verifying a simpler condition: If the random variables $\*x_i$ generating the observations $x_i$ are independent and identically distributed, then it suffices
that Assumption \ref{AS:A_exists} holds. That is to say, it suffices that $\pi$ is \textit{structurally} less informative than $\mathcal{Q}^{\*\theta}$ in the sense of Remark \ref{remark:A_exists}.
%
%
\begin{theorem}[\GVI consistency under independence]
    If Assumptions \ref{AS:min_exists}, \ref{AS:dirac}, \ref{AS:D}, \ref{AS:suitable}, \ref{AS:finite_solution_exists} and \ref{AS:A_exists} hold and $x_i \overset{iid}{\sim} \*x_1$, then the \GVI posteriors are consistent. I.e., $q_n \overset{\mathcal{D}}{\to} \delta_{\*\theta^{\ast}}$ $\mu$-almost surely, where $\mu$ is the probability measure on $\*x_1$.
    \label{thm:indep_consistency_sketch}
\end{theorem}
\begin{proofsketch}
    If $x_i \overset{iid}{\sim} \*x_i$, then it also follows that $\mathbb{E}_{\overline{q}_n}\left[\ell(\*\theta, x_i)\right] \overset{iid}{\sim} \mathbb{E}_{\overline{q}_n}\left[\ell(\*\theta, \*x_1)\right]$. 
    Further, one can show that the triangular array $\{ \{ 
    \mathbb{E}_{\overline{q}_n}\left[\ell(\*\theta, x_i)\right] \}_{i=1}^{n}\}_{n=1}^{\infty}$ satisfies a strong law of large numbers, which is to say that $\frac{1}{n}\sum_{i=1}^n\mathbb{E}_{\overline{q}_n}\left[\ell(\*\theta, x_i)\right] \to \mathbb{E}_{\mu}\left[\ell(\*\theta^{\ast}, \*x)\right]$, $\mu$-almost surely.
    Assumption \ref{AS:A_exists} plays a crucial part in the proof of this law of large numbers. Specifically, it guarantees that a dominated convergence theorem holds outside an (arbitrarily large) compact set $A$ around the optimal value $\*\theta^{\ast}$. 
    As a consequence, it suffices to show that $\{ \{ 
    \mathbb{E}_{\overline{q}_n}\left[\ell(\*\theta, x_i)\cdot1_A(\*\theta)\right] \}_{i=1}^{n}\}_{n=1}^{\infty}$ satisfies a strong law of large numbers and specifically that $\frac{1}{n}\sum_{i=1}^n\mathbb{E}_{\overline{q}_n}\left[\ell(\*\theta, x_i)\cdot1_A(\*\theta)\right] \to \mathbb{E}_{\mu}\left[\ell(\*\theta^{\ast}, \*x)\right]$, $\mu$-almost surely.
    As it turns out, this is straightforward due to the mild regularity conditions of Assumption \ref{AS:min_exists}, which imply that $\ell(\*\theta, x_i)\cdot1_A(\*\theta)$ is upper-bounded in $\*\theta$.
    Together, this entails that $\varepsilon_n$ indeed goes to zero ($\mu$-almost surely), which is equivalent to proving consistency.
\end{proofsketch}

Assumption \ref{AS:A_exists} upon which the last result is based holds for a variety of variational families. For example, Lemma \ref{lemma:A_exists_for_MFN} shows that it holds for virtually all priors of practical interest if $\mathcal{Q}^{\*\theta}$ is the mean field normal family.
Combining this with Theorem \ref{thm:indep_consistency_sketch} then  yields the following result.
\begin{corollary}[\GVI consistency for the mean field normal family]
    If Assumptions \ref{AS:min_exists}, \ref{AS:dirac}, \ref{AS:D}, \ref{AS:suitable}, \ref{AS:finite_solution_exists} hold,  $\pi$ has monotonically decaying tails vanishing at most as fast as Gaussian tails outside a compact set $A$, $x_i \overset{iid}{\sim} \*x_1$ and $\mathcal{Q}^{\*\theta}$ is the mean field normal family, then $q_n \overset{\mathcal{D}}{\to} \delta_{\*\theta^{\ast}}$ $\mu$-almost surely, where $\mu$ is the probability measure on $\*x_1$.
    Priors satisfying this conditions include Gaussians, Student's $t$ distributions, Laplace distributions as well as finite mixtures of such priors.
    \label{corollary:main_corollary_1}
\end{corollary}
\begin{proofsketch}
    %
    First, note that it is clear that the first step in the proof of Theorem \ref{thm:mainres_sketch} holds: $\overline{q}_n \overset{\mathcal{D}}{\toL}\delta_{\*\theta^{\ast}}$.
    Further, the elements in the sequence $\overline{q}_n$ for $\*\Theta \subseteq \mathbb{R}^{D}$ are products of $D$ independent Gaussians: $\overline{q}_n(\*\theta) = \prod_{d=1}^DN(\*\theta_d|\*\mu_{d,n}, \*\sigma_{d,n})$.
    Clearly, since $q_n$ collapses to a point mass, it holds that $(\*\mu_{d,n}, \*\sigma_{d,n}) \to (\*\theta^{\ast}_d, 0)$ as $n\to\infty$.
    To finish the proof, all one needs is a formalization of the immediate intuition that eventually, the tails of $\overline{q}_n$ outside of some compact set $A$ containing $\*\theta^{\ast}$ will have retreated sufficiently to be dominated by $\pi$.
    This holds whenever $\pi$ has monotonically decaying tails outside of $A$ and if these tails decay at most as fast as Gaussian tails outside of $A$. 
    It is immediately clear that this holds for Gaussians, Student's $t$ distributions and Laplace distributions. For $\pi$ being a mixture, it is similarly obvious once one increases $A$ sufficiently so that $\pi$'s tails are monotonically decreasing outside of it.
\end{proofsketch}



%
While the mean field normal family is an interesting special case, one is often interested in better approximations. For variational families that are more expressive than the mean field normal family, Assumption \ref{AS:A_exists} can be hard to verify. 
To elegantly circumvent this issue, Appendix \ref{sec:base_family_stragey} outlines the \textit{base family strategy}: If one can show that a base family $\mathcal{Q}_1^{\*\theta}$ yields consistent posteriors, then a strictly larger variational family $\mathcal{Q}_2^{\*\theta}$ will also yield consistent posteriors. 
Formally, this finding is summarized in the following Theorem.

\begin{theorem}[\GVI consistency via base family strategy]
    Suppose Assumptions \ref{AS:min_exists}, \ref{AS:dirac}, \ref{AS:D}, \ref{AS:suitable}, \ref{AS:finite_solution_exists} hold for the \GVI problem $P(\ell, D, \mathcal{Q}_1)$ with solution sequence $\{q_{n,1}\}_{n=1}^{\infty}$.
    Further, let $\mathcal{Q}_1 \subseteq \mathcal{Q}_2$ and $\{q_{n,2}\}_{n=1}^{\infty}$ be the solution sequence of $P(\ell, D, \mathcal{Q}_2)$.
    %
    %
    If $q_{n,1} \overset{\mathcal{D}}{\to} \delta_{\*\theta^{\ast}}$ $\mu$-almost surely (in the $\mu$-probability limit), then it also holds that $q_{n,2} \overset{\mathcal{D}}{\to} \delta_{\*\theta^{\ast}}$ $\mu$-almost surely (in the $\mu$-probability limit) as $n\to\infty$.
\end{theorem}

\begin{proofsketch}
    Recalling notation, one has $\mathcal{Q}_1 = \mathcal{Q}_1^{\*\theta} \times \mathcal{Q}_{n,1}^z$ and $\mathcal{Q}_2 = \mathcal{Q}_2^{\*\theta} \times \mathcal{Q}_{n,2}^z$ an encompassing family, which is to say that $\mathcal{Q}_1 \subseteq \mathcal{Q}_2$.  Also recall that the second step for the proof in Theorem \ref{thm:mainres_sketch} (detailed in \textbf{Appendix \ref{sec:varepsilon-minimizers}})
    shows that
    \begin{IEEEeqnarray}{rCl}
        \overline{F}_n(q_{n,1}) 
        & \leq &
        \inf_{q\in\mathcal{Q}_1^{\*\theta}}\overline{F}_n(q) + \varepsilon_{n,1}.
        \nonumber
    \end{IEEEeqnarray}
    Further, the third step (detailed in \textbf{Appendix \ref{sec:varepsilon-convergence}})
    shows that $\varepsilon_{n,1} {\toL} 0$ $\mu$-almost surely  (in the $\mu$-probability limit) as $n\to \infty$, which implies that $q_{n,1}\overset{\mathcal{D}}{\to} \delta_{\*\theta^{\ast}}$ holds $\mu$-almost surely (in the $\mu$-probability limit).
    Additionally, one now observes that $\mathcal{Q}_2^{\*\theta}$ is a larger solution space than $\mathcal{Q}_1^{\*\theta}$, so that one gets 
    \begin{IEEEeqnarray}{rCl}
        \overline{F}_n(q_{n,2}) 
        & \leq &
        \inf_{q\in\mathcal{Q}_2^{\*\theta}}\overline{F}_n(q) + \varepsilon_{n,2}
        \nonumber
    \end{IEEEeqnarray}
    for a new error sequence $\varepsilon_{n,2}$. 
    While it is possible that $\varepsilon_{n,2}>\varepsilon_{n,1}$ for some or even all $n$, one can show that $\mu$-almost surely  (in the $\mu$-probability limit), it holds that $|\varepsilon_{n,2} - \varepsilon_{n,1}|\to 0$  as $n\to\infty$. %
    In turn, this shows that $q_{n,2}\overset{\mathcal{D}}{\to} \delta_{\*\theta^{\ast}}$ holds $\mu$-almost surely (in the $\mu$-probability limit).
\end{proofsketch}

This result is not only generic, intuitive and powerful, but also very useful in practice. For example, applying it to the mean field variational family immediately yields the following result.

\begin{corollary}[\GVI consistency for families encompassing the mean field normal]
    If the conditions of Corollary \ref{corollary:main_corollary_1} 
    are satisfied with a variational family $\mathcal{Q}^{\*\theta}$ large enough to encompass the mean field normal family, then the \GVI posteriors are consistent. That is to say, $q_n \overset{\mathcal{D}}{\to} \delta_{\*\theta^{\ast}}$ $\mu$-almost surely, where $\mu$ is the probability measure on $\*x_1$.
    Examples for $\mathcal{Q}^{\*\theta}$ are mixtures of normals, neural networks parameterizing a normal as well as Gaussian Processes.
\end{corollary}

The results presented thus far are interesting, but have relied on independence. In fact, the  proof strategies developed in Appendix \ref{sec:proof_details} are far more widely applicable and do \textbf{not} depend on the independence assumption per se.
In particular and as emphasized in the proof of Theorem \ref{thm:mainres_sketch}, all that is needed is that a law of large numbers can be shown to hold for the triangular array $\left\{ \left\{ \mathbb{E}_{\overline{q}_n}\left[ 
        \ell(\*\theta, x_i)
        \right] \right\}_{i=1}^n \right\}_{n=1}^{\infty}$.
While this is especially easy under independence, there are many situations under which laws of large numbers hold even under severe dependence. 
%
%
To show the practical relevance of this observation, the following Theorem collects some of the more easily verifiable conditions that guarantee consistency under dependence. 
All conditions listed are simplified special cases of the ones presented in Corollaries \ref{corollary:summary_consistency} and \ref{corollary:summary_consistency_dependence}.

\begin{theorem}[\GVI consistency under dependence]
    Suppose that for a given \GVI problem $P(\ell, D, \mathcal{Q})$ with potentially dependent observations $x_{1:n}$, Assumptions \ref{AS:min_exists}, \ref{AS:dirac}, \ref{AS:D}, \ref{AS:suitable} and \ref{AS:finite_solution_exists} hold. Suppose that additionally, one of the following holds: 
    \begin{itemize}
        \item[(a)] $\ell(\*\theta, \*x) \leq M$ for some constant $M$, $\mu$-almost surely;
        \item[(b)] $\ell(\*\theta, \*x)$ is jointly continuous in $\*\theta$ and $\*x$, $\mu$-almost surely and both $\mathcal{X}$ and $\*\Theta$ are compact;
        \item[(c)] $\frac{1}{n}\sum_{i=1}^n\ell(\*\theta, x_i)$ converges to $\mathbb{E}_{\mu}\left[\ell(\*\theta, \*x)\right]$ both $\mu$-almost surely and uniformly over $\*\Theta$.
    \end{itemize}
    Then, the \GVI posteriors are strongly consistent. That is to say, $q_n \overset{D}{\to} \delta_{\*\theta^{\ast}}$ $\mu$-almost surely.
    \label{thm:consistency_dependence_sketch}
\end{theorem}
\begin{proofsketch}
    Conditions (a), (b) and (c) are all variations on the same theme: All of them provide different notions of boundedness. 
    While (a) and (b) impose a form of boundedness directly upon $\ell$, condition (c) is more subtle and imposes a stochastic form of boundedness on the difference $\sup_{\*\theta \in \*\Theta}|\frac{1}{n}\sum_{i=1}^n\ell(\*\theta, x_i) - \mathbb{E}_{\mu}\left[\ell(\*\theta, \*x)\right]|$.
    Either way, imposing boundedness is a way to enable application of a dominated convergence theorem.
    In turn, this enables proving that $\mu$-almost surely, $\frac{1}{n}\sum_{i=1}^n\mathbb{E}_{\overline{q}_n}\left[ 
        \ell(\*\theta, x_i)
        \right] \to \mathbb{E}_{\mu}\left[ \ell(\*\theta^{\ast}, \*x)\right]$, which is to say that $\varepsilon_n \to 0$.
\end{proofsketch}



\section{Experiments}

Having formally proved \GVI consistency, the next section provides a short simulation study to verify the findings.
First, convergence properties are studied for a correctly specified Bayesian Linear Regression (\BLR) model using a range of divergences $D$ to penalize deviations of the posterior from the prior belief. 
Second, two scenarios are studied for a Bayesian Mixture Model (\BMM). In the first scenario, a well-specified is contrasted with an ill-specified prior belief. In the second scenario, the effect of model misspecification is investigated. 
For all examples, the variational family within which the posteriors are computed is a dimension-wise independent normal distribution (the so-called \textit{mean-field variational family}) which for $p_{\mathcal{N}}(\*\theta_j|m_j, s_j)$ the density of a normal distribution on $\*\theta_j$ with mean $m_j$ and standard deviation $s_j$ is given by
\begin{IEEEeqnarray}{rCl}
    \mathcal{Q}_{\text{MFN}} = \left\{ \prod_{j=1}^d p_{\mathcal{N}}(\*\theta_j|m_j, s_j): m_j \in \mathbb{R}, s_j \in \mathbb{R}_{>0}  \right\}.
    \nonumber
\end{IEEEeqnarray}

\subsection{Bayesian Linear Regression (\BLR)}

\begin{figure}[t!]
 \begin{center}
 \includegraphics[trim= {3.0cm 0.0cm 3.9cm 0.0cm}, clip,   width=1\columnwidth]{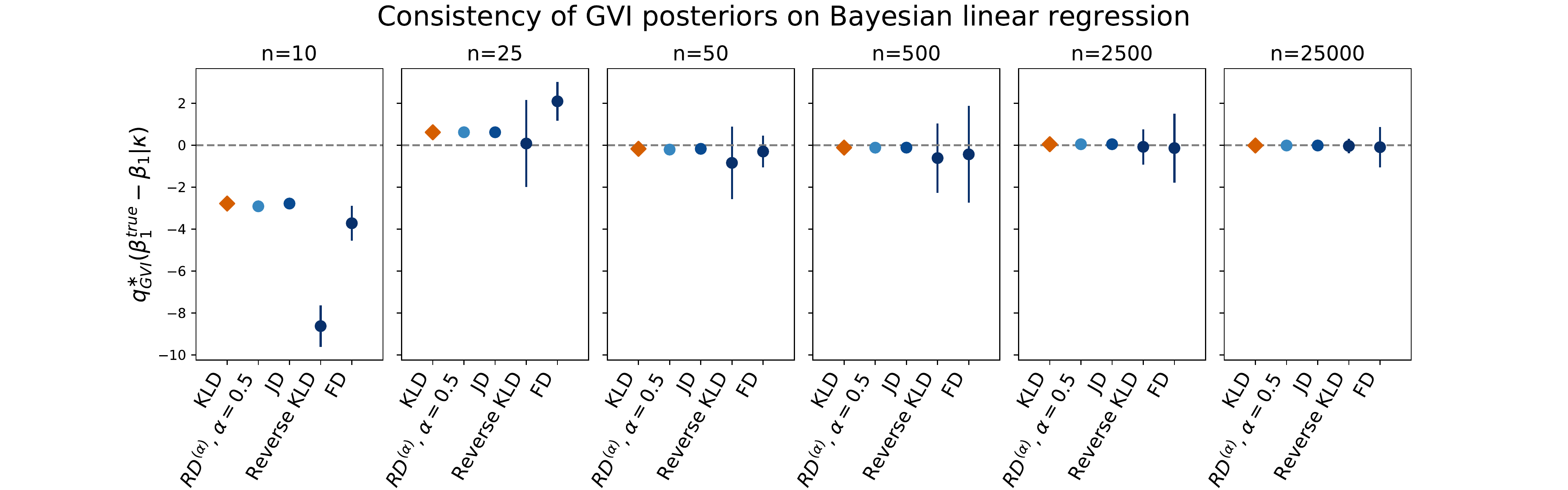}
 \caption{
 Depicted are the inferred \VIColor{\textbf{\VI}} and  \GVIColor{\textbf{\GVI}} posteriors for $\*\beta_1$ using the uncertainty quantifiers given in Table \ref{table:divergences}. 
Because all inferred posterior beliefs are normals, dots are used to mark out the posterior mean and whiskers to denote the posterior standard deviation.
All posteriors are re-centered around the true value of $\*\beta_1$, so that the $y$-axis shows how far the posterior belief is from the truth.
 }
 \label{Fig:BLR_consistency}
 \end{center}
 \end{figure}

For all examples of the \BLR model, the $i$-th observation is given by $x_i^0 = (y_i, \widetilde{x}_i)$ and generated as
\begin{IEEEeqnarray}{rCl}
	y_i 
	& = & 
	\widetilde{x}_i^T \*\beta + \varepsilon_i 
	\nonumber \\
	\varepsilon_i 
	& \overset{i.i.d.}{\sim} & 
	\mathcal{N}\left(0, \sigma^2 \right).
	\nonumber
\end{IEEEeqnarray}
The regressors $\widetilde{x}_i\in\mathbb{R}^{20}$ themselves are drawn from independent normal distributions. The variance parameter is set to be  $\sigma^2 = 25$ and  the coefficients $\*\beta\in\mathbb{R}^{20}$ are fixed across all simulations and given in the Appendix. 
Much like the regressors $\widetilde{x}_i$, the values of $\*\beta$ were generated as independent draws from a normal with mean $3$ and standard deviation $10$ and are all contained in the interval $(-12.46, 17.85)$. 
Across all simulations, the loss is the accurately specified negative log likelihood of the \BLR model. For the $i$-th observation, $\*\theta = (\*\beta, \sigma^2)$ and $p_{\mathcal{N}}(\mu, \sigma^2)$ the likelihood of a normal with mean $\mu$ and variance $\sigma^2$, the loss is thus given by 
\begin{IEEEeqnarray}{rCl}
	\ell(\*\theta, x_i^o) 
	& = & 
	- \log p_{\mathcal{N}}(y_i - \widetilde{x}_i^T \*\beta, \sigma^2).
	\nonumber
\end{IEEEeqnarray}
Similarly, the prior is fixed across all settings to be the fully factorized normal distribution 
\begin{IEEEeqnarray}{rCl}
	\pi(\*\theta) & = &
	p_{\mathcal{N}}(\log(\sigma^2), 10^2) \prod_{d=1}^{20} p_{\mathcal{N}}(\*\beta_d, 10^2).
	\nonumber
\end{IEEEeqnarray}
In contrast to the loss and the prior, the uncertainty quantifier $D$ is varied. Table \ref{table:divergences} shows all uncertainty quantifiers used in the simulations. 
The Appendix explains their closed forms if both $\pi$ and $q$ are fully factorized normals. 
Varying the total number of observations $n$, Figure \ref{Fig:BLR_consistency} compares \GVI posteriors based on these divergences to demonstrate that frequentist consistency is satisfied by any of the resulting \GVI posteriors.

The plot also shows that changing $D$ has considerable impact on both uncertainty quantification as well as convergence speed. 
%
%
For example, the $\alpha$-divergence and R\'enyi's $\alpha$-divergence (\RAD) produce posterior beliefs that converge virtually at the same speed as the standard \VI posterior.
A very different behaviour is exhibited by the reverse \KLD or the Fisher Divergence (\FD): These choices produce posterior beliefs that encode a strictly larger degree of uncertainty than the \KLD.  
%
%
This can be interpreted to mean that they encode a much stronger trust in the prior than the \KLD.

\begin{table}[th!]
\begin{center}
\begin{small}
\begin{tabular}{ p{4.5cm}p{8.4cm}}
    \multicolumn{1}{l}{Name \& abbreviation} & 
    \multicolumn{1}{l}{Definition of $D(q\|\pi)$}  
	 \\[0.1cm]
    \hline\hline &\\[-0.15cm]
    		Kullback-Leibler Divergence (\KLD) &
    		$\KLD(q\|\pi) = \mathbb{E}_{q}\left[ 
    			\log\left(\dfrac{q(\*\theta)}{\pi(\*\theta)}\right)	
    		\right]$ 
    		 \\[0.35cm]
    		 Reversed \KLD &
    		$\KLD(\pi\|q) = \mathbb{E}_{\pi}\left[ 
    			\log\left(\dfrac{\pi(\*\theta)}{q(\*\theta)}\right)	
    		\right]$ 
    		 \\[0.35cm]
    		R\'enyi's $\alpha$-divergence (\RAD) &
    		$\RAD(q\|\pi) = \dfrac{1}{\alpha-1}\log\left( 
    			\mathbb{E}_q\left[ 
    				 \left(
    				 	\dfrac{q(\*\theta)}{\pi(\*\theta)}
    				  \right)^{\alpha - 1} \right] \right)$, 
    				  $\alpha> 0, \alpha \neq 1$
    		\\[0.5cm]
    		$\alpha$-divergence (\AD) &
    		$\AD(q\|\pi) = \dfrac{1}{\alpha(\alpha-1)}
    			\mathbb{E}_q\left[ 
    				 \left(
    				 	\dfrac{q(\*\theta)}{\pi(\*\theta)}
    				  \right)^{\alpha - 1} - 1\right]$, 
    				 $\alpha> 0, \alpha \neq 1$
    		\\[0.5cm]
    		Jeffrey's Divergence (\JD) &
    		$\JD(q\|\pi) = \KLD(q\|\pi) + \KLD(\pi\|q)$
    		\\[0.15cm]
    		%
    		Fisher's Divergence (\FD) &
    		$\FD(q\|\pi) = \mathbb{E}_{q}\left[ \|
    			\nabla_{\*\theta}\log q(\*\theta) - \nabla_{\*\theta} \log \pi(\*\theta) \|_2^2	
    		\right]$ 
    		\\[0.35cm]
    \hline \hline  \\
\end{tabular}
\end{small}
{\renewcommand{\arraystretch}{1.2}
\caption{A range of divergences satisfying Assumption \ref{AS:D}. Figure \ref{Fig:BLR_consistency} demonstrates that the choice of $D$ crucially impacts the speed at which the \GVI posterior converges to a point mass at the population-optimal value.
}
\label{table:divergences}
}
\end{center}
\vskip -0.1in
\end{table}

\subsection{Bayesian Mixture Model (\BMM)}

Though the consistency results also hold in the presence of latent variables, the \BLR example of the previous section only depended on observables $x^o_i$.
The next section extends this to Bayesian Mixture Models (\BMM{}s), which are latent variable models.
Before empirically demonstrating consistency however, the section first introduces the experimental setting and demonstrates the potential appeal for using $D\neq \KLD$. 

\subsubsection{Experimental setup}

Throughout, $n$ observations are generated from the $d$-dimensional \BMM with two equally likely normal mixture components $z=1,2$ with dimension-wise unit variance and mean given  by
\begin{IEEEeqnarray}{rCl}
    \*\mu^z  = (\mu_1^z, \mu_2^z, \dots \mu_d^z)^T 
    & = & 
    \begin{cases}
        2 \cdot e_d & \text{ if } z=1 \\
        -2 \cdot e_d & \text{ if } z=2
    \end{cases},
    \nonumber
\end{IEEEeqnarray}
where $e_d = (1, 1, \dots 1)^T$ is the $d$-dimensional column vector of ones. 
The $n$ observations $x_i^o$ are drawn with equal probability from two mixture, meaning that
\begin{IEEEeqnarray}{rCl}
    \*z_i & \overset{i.i.d.}{\sim} & \text{Bernoulli}(0.5) \nonumber \\
    x_i^o|\{\*z_i = z_i\} & \overset{i.i.d.}{\sim}& \mathcal{N}(x_i^o|\*\mu^{z_i}, \*I_d).
    \label{eq:BMM_data_generation}
\end{IEEEeqnarray}
Notice in particular that this generates $n$ latent variables $\*z_{1:n}$ that indicate mixture memberships for $x_{1:n}^o$, but are unobserved.
With this, inference is conducted on $\*\mu^c$ for $c=1,2$ via the negative log likelihood loss of the correct model.
For $\*\theta = (\*\mu^1, \*\mu^2)$, this is given by
\begin{IEEEeqnarray}{rCl}
    \ell(\*\theta, x_i^o, \*z_i) & = &
    - \log p_{\mathcal{N}}(x_i^o|\*\mu^{\*z_i}, \*I_d).
    \nonumber
\end{IEEEeqnarray}
Two variations of this experiment are considered: In Section \ref{sec:experiment_robust_D}, the benefits of alternative choices of $D$ are explored for the fixed number of observations $n=50$.
To this end, $B=100$ artificial data sets are generated according to the above description. 
The second variation of the experiment is given in Section \ref{sec:experiment_consistency_across_D} and explores frequentist consistency as $n\to\infty$. On top of that, it also studies the effects that model misspecification and robust scoring rules have on inference outcomes.
In particular, two settings are explored: In the first case, the data is generated as above. In the second case, additional noise is injected. Specifically, after $x_i^o$ is generated according to eq. \eqref{eq:BMM_data_generation}, inference is based on the possible polluted observation $\widetilde{x}_i^o$ generated as
\begin{IEEEeqnarray}{rCl}
    \widetilde{x}_i^o & = & x_i^o + u_i\cdot\eta_i\cdot e_d \nonumber \\
    u_i & \overset{i.i.d.}{\sim} & \text{Bernoulli}(0.05) \nonumber \\
    \eta_i & \overset{i.i.d.}{\sim} & \mathcal{N}(10, \sqrt{3}). \label{eq:BMM_data_generation_contamination}
\end{IEEEeqnarray}

\subsubsection{Why use $D\neq \KLD$?}
\label{sec:experiment_robust_D}

\begin{figure}[t!]
 \begin{center}
 \includegraphics[trim= {1.2cm 3.5cm 3.0cm 0.0cm}, clip,   width=1\columnwidth]{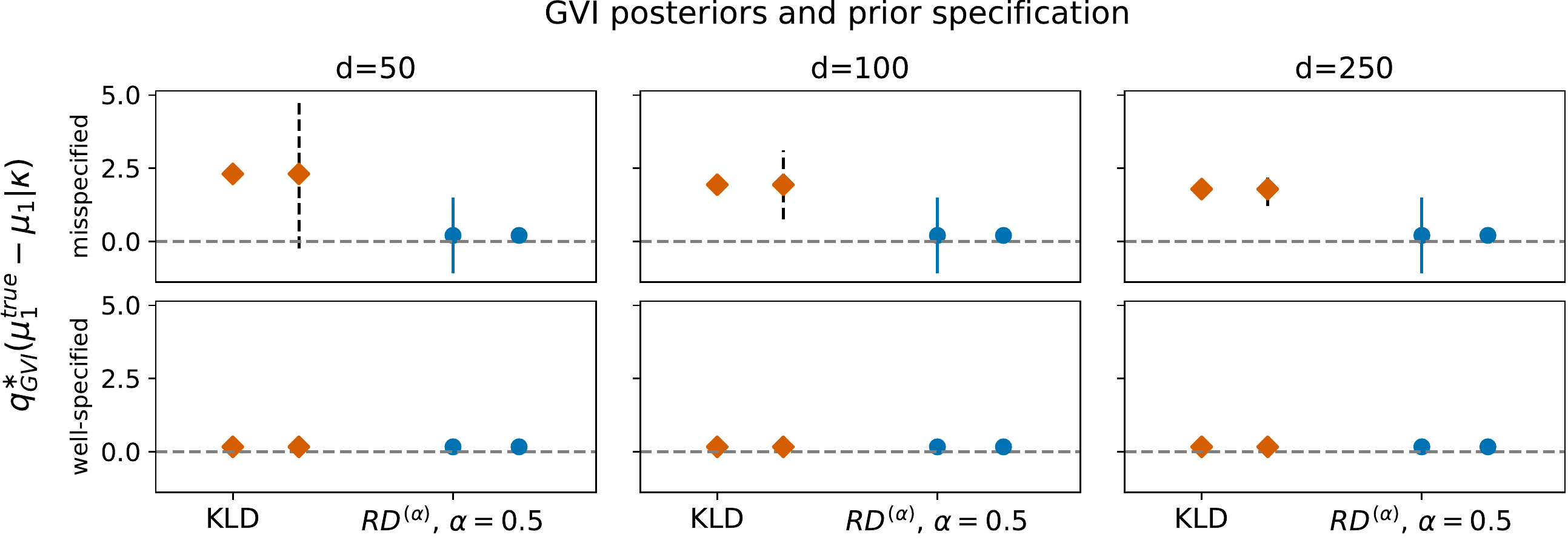}
 \caption{
 The \textbf{first column} of each setting depicts the inferred \VIColor{\textbf{\VI}} and  \GVIColor{\textbf{\GVI}} posteriors for $\*\theta$ in the \BMM{} of eq. \eqref{eq:BMM_data_generation}.
 Here, the \GVIColor{\textbf{\GVI}} posteriors use $D=\RAD$ for $\alpha = 0.5$. 
 All inferred posterior beliefs are normals, so dots and whiskers mark posterior means and standard deviations.
 The posterior are re-centered so that the $y$-axis measures the magnitude by which the posterior belief deviates from the truth.
 The \textbf{second column} of each setting shows the inferred posterior mean and its standard error across the 100 data sets on which the experiment was run. The plots clearly show that the adverse effect of the prior stabilizes as the number $d$ of affected parameters increases.
 }
 \label{Fig:BMM_why_use_GVI}
 \end{center}
 \end{figure}

If the prior is poorly specified, $D=\KLD$ will produce posterior beliefs that place the same weight on the prior as they do on the data. 
In contrast, robust alternatives to the \KLD do not suffer this problem: They can produce posterior beliefs that take the prior into account, but are robust to prior misspecification, see also \citet{GVI}. 
To illustrate the phenomenon empirically, the next experiment compares the \KLD with R\'enyi's $\alpha$-divergence (\RAD) for $\alpha=0.5$ under two settings: A well-specified prior $\pi_1(\*\theta)$ and a misspecified prior $\pi_2(\*\theta)$, which are given by
\begin{IEEEeqnarray}{rCl}
    \pi_1(\*\theta) & = &
    p_{\mathcal{N}}\left(\*\theta|0_d, \sqrt{10}\*I_d\right)
    \nonumber \\
    \pi_2(\*\theta) & = &
    p_{\mathcal{N}}\left(\*\theta|-10 \cdot e_d, \sqrt{0.1}\*I_d\right)
    \nonumber
\end{IEEEeqnarray}
In a nutshell, \RAD is a robust alternative for $\alpha \in (0,1)$ \citep{ABCdiv}, but recovers the \KLD for $\alpha \to 1$. This gives the magnitude of $\alpha$ an interpretation as the degree of prior robustness.
To evaluate the experiments, 100 data sets are generated with $n=50$ observations each. Across these, Figure \ref{Fig:BMM_why_use_GVI} reports the average posterior computed as
\begin{IEEEeqnarray}{rCl}
    \mathcal{N}\left(\bar{m}, \bar{s}\right), \quad \bar{m} = \frac{1}{100}\sum_{j=1}^{2d}\sum_{b=1}^Bm_{b,j}, \quad  \bar{s} = \frac{1}{100}\sum_{j=1}^{2d}\sum_{b=1}^Bs_{b,j}. \nonumber
\end{IEEEeqnarray}
Here, $s_{b,j}$ corresponds to the standard deviation computed for the $j$-th dimension of the mean field normal posterior on the $b$-th artificial data sets.
Similarly, $m_{b,j}$ corresponds to the mean of the same parameter posterior, albeit re-centered around the true value of the inferred parameter.  

As Figure \ref{Fig:BMM_why_use_GVI} shows, \RAD is an interesting alternative to the \KLD in finite samples: If the prior is misspecified (top row), the \KLD produces belief distributions that take the prior too strongly into account and are far from the truth. 
In contrast, the \RAD provides both prior robustness as well as better uncertainty quantification under misspecification.
At the same time,  \RAD has no tangible disadvantage relative to the \KLD if the prior is well-specified (bottom row).

 \begin{figure}[hp!]
 \vskip -1.0cm
 \begin{center}
 \includegraphics[trim= {0.5cm 0.0cm 2.5cm 0.0cm}, clip,   width=0.925\columnwidth]{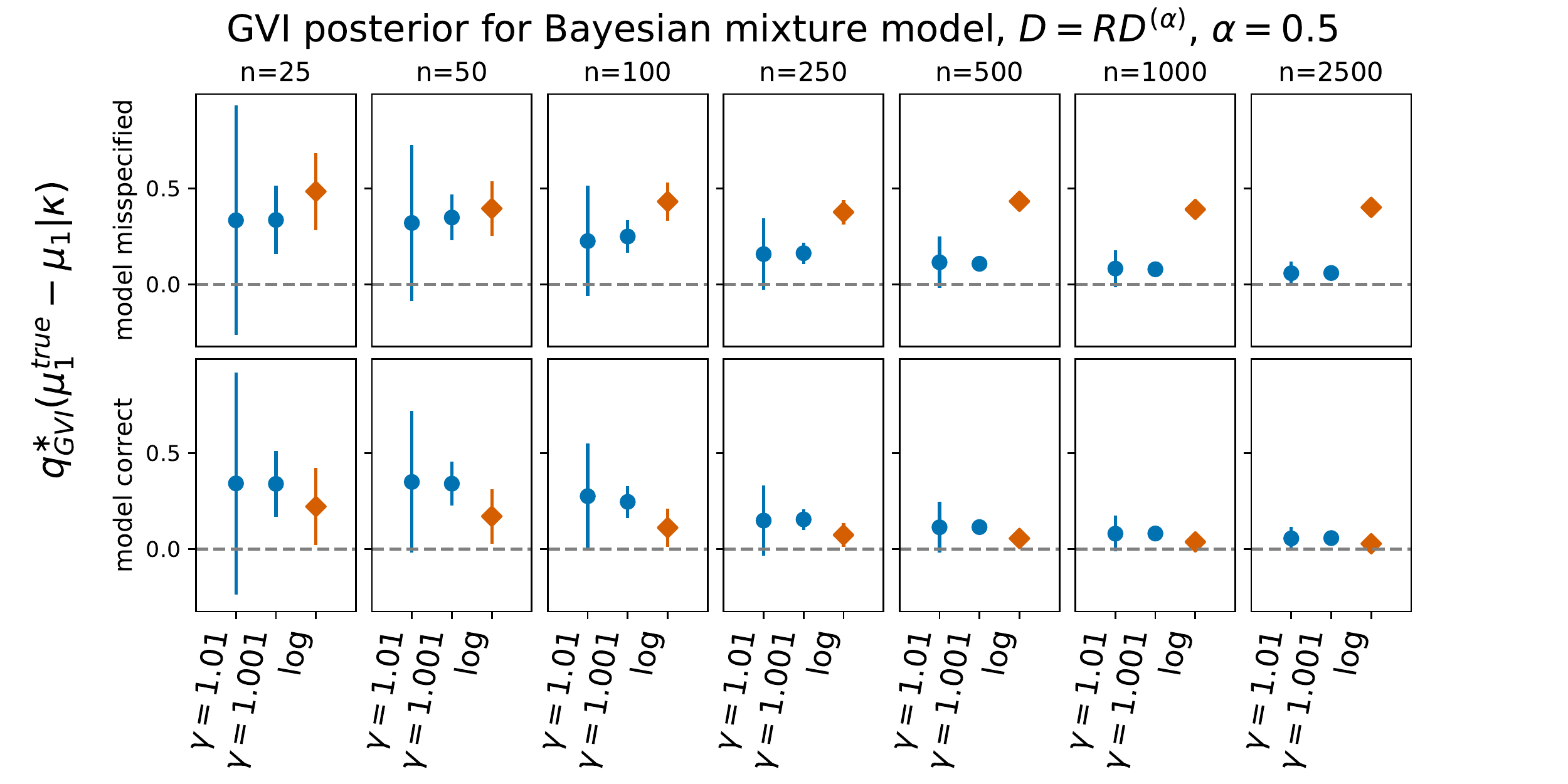}
 \includegraphics[trim= {0.5cm 0.0cm 2.5cm 0.0cm}, clip,   width=0.925\columnwidth]{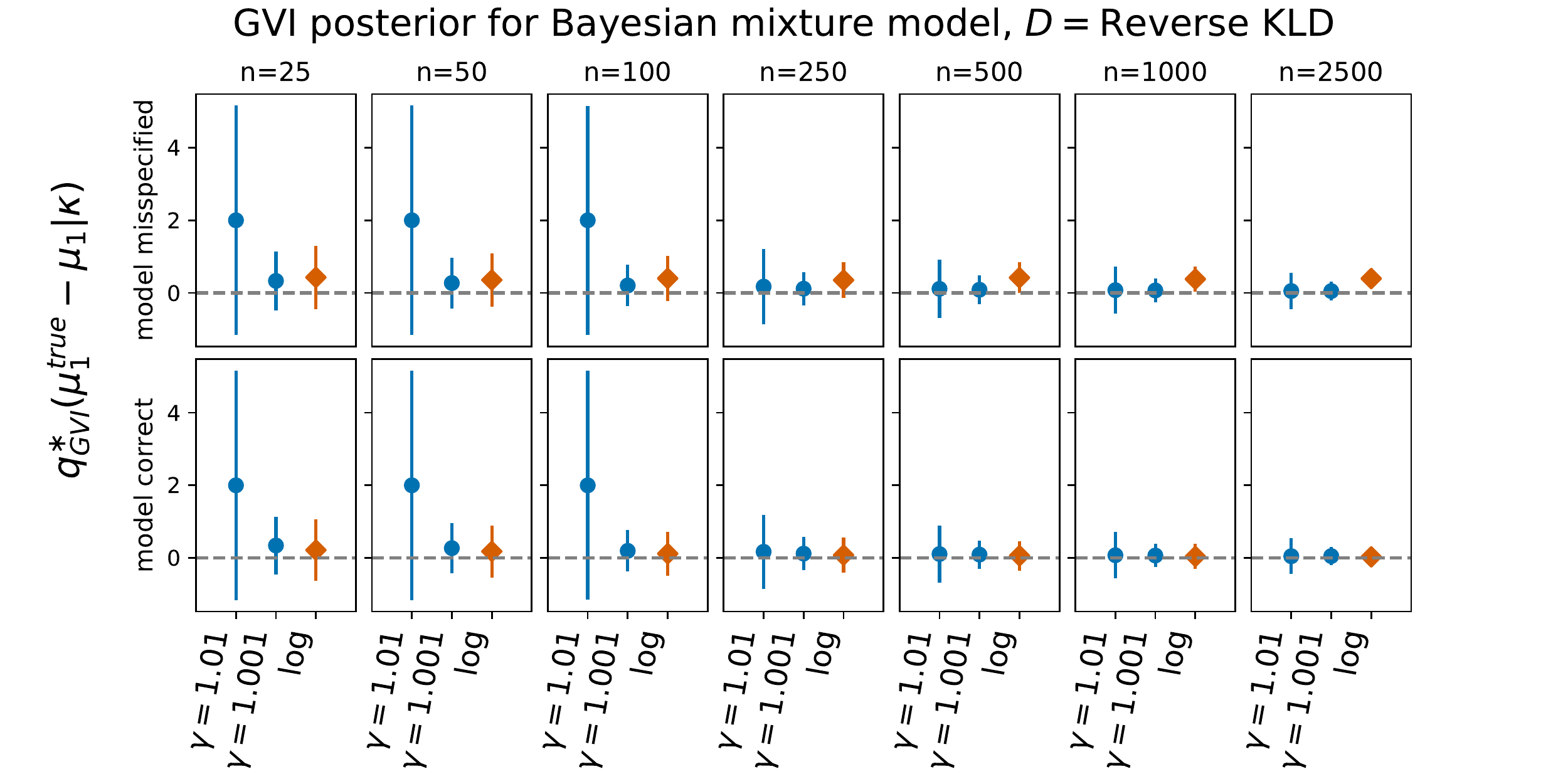}
 \includegraphics[trim= {0.5cm 0.0cm 2.5cm 0.0cm}, clip,   width=0.925\columnwidth]{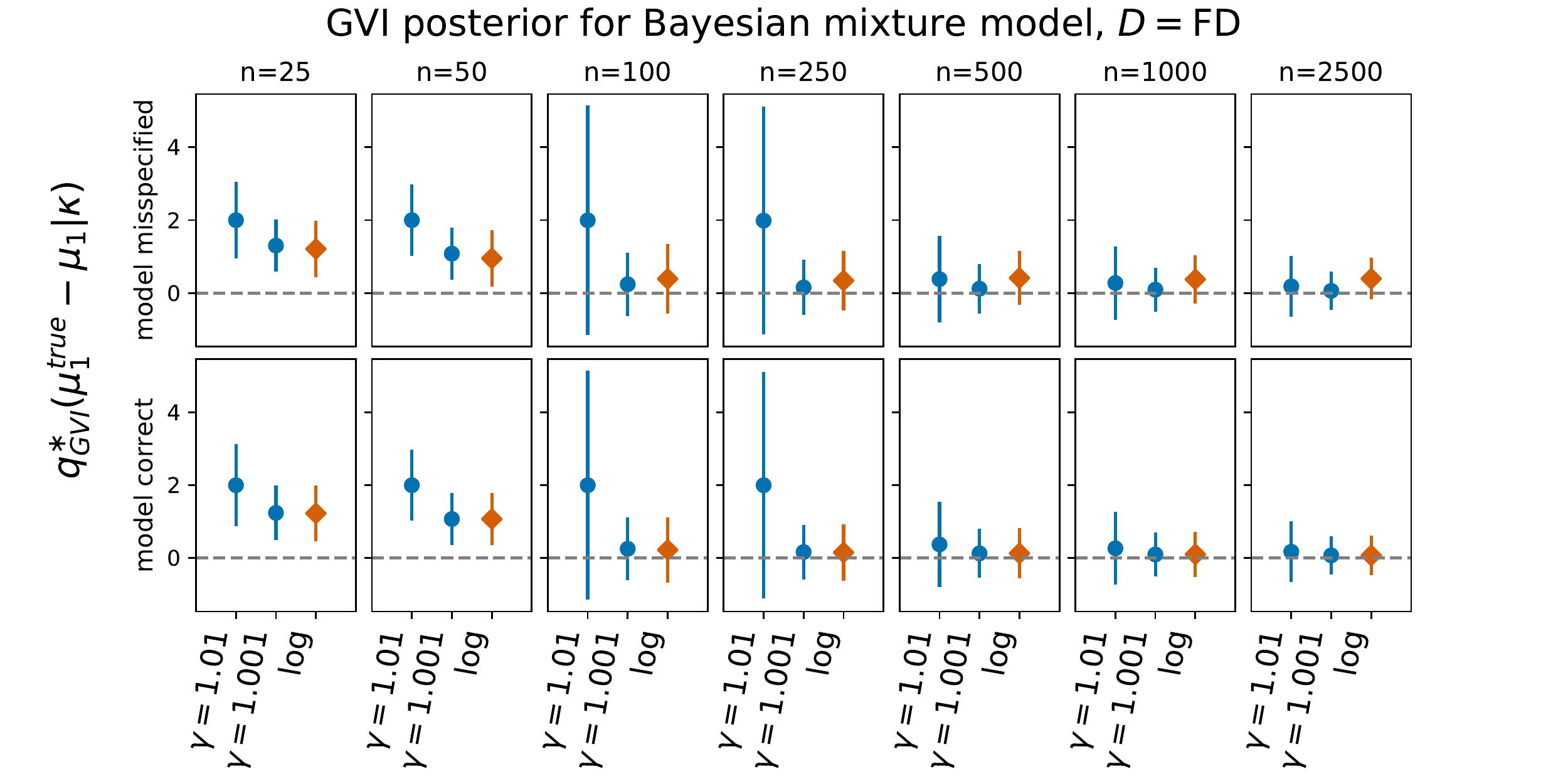}
 \caption{
 Depicted are the inferred \VIColor{\textbf{\VI}} and  \GVIColor{\textbf{\GVI}} posteriors for $\mu$. Here, the \GVIColor{\textbf{\GVI}} posteriors use $D=\RAD$ for $\alpha = 0.5$. 
Because all inferred posterior beliefs are normals, dots are used to mark out the posterior mean and whiskers to denote the posterior standard deviation.
All posteriors are re-centered around the true value of $\*\beta_1$.
 }
 \label{Fig:BMM_consistency}
 \end{center}
 \end{figure}

 \subsubsection{Consistency}
 \label{sec:experiment_consistency_across_D}
 
 Having motivated the use of $D\neq \KLD$ in finite samples, the next experiment investigates convergences speed as a function of the choice for $D$.
 %
 We use the opportunity to compare consistency in two settings: If the model is correctly specified as in eq. \eqref{eq:BMM_data_generation} and when it is misspecified via eq. \eqref{eq:BMM_data_generation_contamination}.
 %
In this context of misspecification, two different types of loss functions are compared: Firstly, the standard negative log likelihood given by
 \begin{IEEEeqnarray}{rCl}
     \ell(\*\theta, x_i) & = & - \log p(x_i|\*\theta). \nonumber
 \end{IEEEeqnarray}
 as well as a robust scoring rule derived from the $\gamma$-divergence \citep{GammaDivSummable} given by
 \begin{IEEEeqnarray}{rCl}
    \Lg(\*\theta, x_i) & = &
        -
			\frac{1}{\gamma-1}p(x_i|\*\theta)^{\gamma-1} \cdot
			  \frac{\gamma}{I_{p, \gamma}(\*\theta)^{-\frac{\gamma-1}{\gamma}}}
			.
		\nonumber
 \end{IEEEeqnarray}
For the relevant background on minimium scoring rule inference, see \citet{dawid2016minimum} and \citet{Jewson}.
It is well-known that the log score is not robust to misspecification \citep[see e.g.][and references therein]{Jewson}. In contrast, $\Lg$ defines a scoring rule that is strongly robust to contamination \citep{GammaDivNotSummable, GammaDivSummable, RobustBayesGamma}. The degree of robustness  is regulated by $\gamma$: While $\gamma > 1$ produces more robust inferences than the $\log$ score, $\Lg$ recovers the $\log$ score as $\gamma \to 1$.
Consequently, one should expect $\Lg$ for $\gamma = 1 + \varepsilon$ for very small values of $\varepsilon>0$ to produce desirable inferences
On the one hand, inferences are nearly as data-efficient as under the $\log$ score if the model is correctly specified. 
On the other hand---and unlike with the log score---the inferences remain reliable under misspecification.

%
Figure \ref{Fig:BMM_consistency} depicts this behaviour and connects it to the consistency findings in the current paper. 
The plot demonstrates two phenomena: 
Firstly, the exact path and speed of the convergence for \GVI posteriors depends on the choice for $D$, especially for small sample sizes. 
Secondly, the overall patterns are the same across all choices of $D$ and are dictated by the choice of $\ell$.
This should not come as a surprise: For $n\to\infty$, the \GVI posteriors concentrate around $\*\theta^{\ast} = \argmin_{\*\theta \in \*\Theta}\mathbb{E}_{\mu}\left[\ell(\*x, \*\theta)\right]$, which does not depend on $D$.
Thus and as the theory predicts, \GVI posteriors collapse to different point masses: While $\Lg$ recovers the true parameter values of eq. \eqref{eq:BMM_data_generation} for $n\to \infty$ in both the misspecified and well-specified setting, the $\log$ only manages to recover the true parameter values in the well-specified setting. In particular,  \GVI posteriors based on the $\log$ score concentrate around a sub-optimal parameter value in the misspecified setting, regardless of the choice for $D$.
%


\section{Discussion}

The current paper gives the first generically applicable consistency result for posterior beliefs generated with Generalized Variational Inference (\GVI).
The results show that virtually all \GVI posteriors of interest contract around the population-optimal parameter value as the number of observations goes to infinity.
In technical terms, the current paper overcomes the substantial challenge of proving consistency for posteriors that do not have any discrepancy-based interpretation as approximations to the exact Bayesian posterior.
This problem is overcome by using an auxiliary objective and techniques from variational calculus together with desirable properties of $\Gamma$-convergent functions.
While results on concentration speed cannot be derived based on this approach, the consistency results are widely applicable for a large range of modularly composable \GVI objectives.
In the future, an interesting question could be a more detailed analysis of \GVI objectives as a function of the uncertainty quantifier $D$. 
Specifically, it is likely that different choices for $D$ yield different concentration rates. 
This would have important ramifications for designing appropriate \GVI objectives: Choosing $D$ based on whether it yields faster (or slower) contraction to the population-optimal value encodes a trade-off between prior-robustness and model-robustness.







\section*{Acknowledgements}

I would like to cordially thank Yulong Yu for fruitful discussions and remarks regarding functional analysis and $\Gamma$-convergence that proved very helpful for the purposes of this paper. 
Similarly, I would like to express my gratitude toward Giuseppe Di Benedetto for important references . 
Further, the work on this paper began to take shape during a three month stay at Duke University under supervision of David Dunson, whom I would like to wholeheartedly thank for getting me interested in questions of Bayesian consistency and giving me useful advice throughout the process.
Beyond that, I would like to thank my PhD supervisor Theodoros Damoulas for giving me the freedom and time to pursue this project.
This work was funded by \EPSRC grant EP/L016710/1 as part of the Oxford-Warwick Statistics Programme (\OxWaSP) as well as by the Facebook Fellowship Programme. 
The work was also supported by the Lloyds Register Foundation programme on Data Centric Engineering through the London Air Quality project and by The Alan Turing Institute for Data Science and AI under \EPSRC grant EP/N510129/1 in collaboration with the Greater London Authority.

\bibliography{library}
\bibliographystyle{apalike}

\newpage
\appendix

\section{Details \& Proofs}
\label{sec:proof_details}

Each of the following sections establishes a key step in the proof of \GVI's consistency. 
First, the auxiliary objectives $\overline{F}_n$ are studied in Section \ref{sec:gamma_convergence}. Specifically, their $\Gamma$-convergence and equi-coerciveness are proved formally, which allows the conclusion that their solution sequences converge. 
Second, Section \ref{sec:varepsilon-minimizers} couples the solutions $q_n$ of the actual objective $F_n$ to the auxiliary objective $\overline{F}_n$ by proving that the former are $\varepsilon_n$-minimizers of the latter.
Thirdly, sufficient conditions are derived in Section \ref{sec:varepsilon-convergence} to guarantee that $\varepsilon_n$ goes to zero ($\mu$-almost surely or in the $\mu$-probability limit) as $n\to\infty$. This finally suffices to show that as desired, the \GVI posteriors collapse to a point mass at the $\mu$-population-optimal value $\*\theta^{\ast}$.
Lastly, a convenient strategy is derived in Section \ref{sec:base_family_stragey} to show that if the \GVI posteriors are consistent when the variational family $\mathcal{Q}_1$ is used, then consistency also holds with any richer family $\mathcal{Q}_2$ encompassing $\mathcal{Q}_1$.

\subsection{Analysis of $\overline{F}_n$}
\label{sec:gamma_convergence}

This section establishes the properties of $\overline{F}_n$ that guarantee that the sequence of its minimizers $\overline{q}_n$ converges to a dirac delta at the $\mu$-population optimal value $\*\theta^{\ast}$  of $\*\theta$.
Lemmas \ref{lemma:coercive_Psi},\ref{lemma:equi_coercive_Fbar} and \ref{lemma:gamma_conv_Fbar} are technical results needed to prove properties about $\overline{F}_n$. Specifically, it is shown that $\overline{F}_n$ is equi-coercive and $\Gamma$-converges to $\mathbb{E}_{q}\left[ \mathbb{E}_{\mu}\left[\ell(\*\theta, \*x)\right]\right]$.
Both findings are pre-requirements needed to establish Lemma \ref{corollary:consistency_Fbar}, which shows that the solutionss oof the auxiliary objective are consistent. In other words, this section proves that so long as Assumptions \ref{AS:min_exists}, \ref{AS:dirac}, \ref{AS:D} and \ref{AS:suitable} hold, $\overline{q}_n \overset{\mathcal{D}}{\toL} \delta_{\*\theta^{\ast}}$.

The next Lemma is a technical result needed to prove equi-coerciveness of $\{\overline{F}_n\}_{n=1}^{\infty}$, which in turn is a pre-requirement to establish consistency of $\overline{q}_n$.

\begin{lemma}
    If Assumption \ref{AS:min_exists} holds,
    $\mathbb{E}_{q}\left[
        \mathbb{E}_{\mu}\left[
            \ell(\*\theta, \*x)
        \right]
    \right]$ is coercive.
    \label{lemma:coercive_Psi}
\end{lemma}
\begin{proof}
    A function is coercive if its sub-level sets are both closed and compact. The sub-level sets of $\Psi(q) = \mathbb{E}_{q}\left[
        \mathbb{E}_{\mu}\left[
            \ell(\*\theta, \*x)
        \right]
    \right]$ are indeed closed since $\Psi(q)$ is lower semi-continuous (see Proposition 1.7 of \citep[][]{gammaConvergence}).
    They are also compact if $\*\Theta$ is compact or if $\mathbb{E}_{\mu}\left[\ell(\*\theta, \*x)\right]$ is coercive in $\*\theta$, as is shown next.
    First, suppose that $(\*\Theta, d)$ is a compact metric space. It is well-known that this immediately implies that $\mathcal{P}(\*\Theta)$ is tight. 
    By Prokhorov's Theorem, this implies that $\mathcal{P}(\*\Theta)$ is compact, which means that any subset of $\mathcal{P}(\*\Theta)$ is compact. Thus, $\Psi(q)$ is coercive whenever $(\*\Theta, d)$ is a compact metric space.
    Second, suppose that $\mathbb{E}_{\mu}\left[\ell(\*\theta, \*x)\right]$ is coercive in $\*\theta$. %
    In other words, $\mathbb{E}_{\mu}\left[\ell(\*\theta, \*x)\right]\to\infty$ as $\|\*\theta\| \to \infty$.
    First, define the sub-level sets as $\mathcal{S}_t = \{q\in\mathcal{Q}^{\*\theta}:\Psi(q)\leq t\}$.
    Since $\mathbb{E}_{\mu}\left[\ell(\*\theta, \*x)\right]$ is coercive in $\*\theta$, for any constant $C\in \mathbb{R}$, there exists ${\*\theta_{C}}$ so that for a sufficiently large ball $B_{\*\theta_{C}}(r_C) = \{\*\theta \in \*\Theta: d(\*\theta, \*\theta_{C}) \leq r_C \}$ of radius $r_C$ around $\*\theta_{C}$, $\mathbb{E}_{\mu}\left[\ell(\*\theta, \*x)\right]\geq C$ for all $\*\theta \notin B_{\*\theta_{C}}(r_C)$. 
    Thus, for any $q \in \mathcal{S}_t$ and for any arbitrarily small and fixed $\Delta > 0$,
    \begin{IEEEeqnarray}{rCCCCl}
        C\int_{\*\Theta \setminus B_{\*\theta_{C}}(r_C)}
        q(\*\theta)d\*\theta
        &\leq &
        \int_{\*\Theta \setminus B_{\*\theta_{C}}(r_C)}
        \mathbb{E}_{\mu}\left[\ell(\*\theta, \*x)\right]q(\*\theta)d\*\theta
        &< &
        t + \Delta <\infty.
        \nonumber
    \end{IEEEeqnarray}
    Rearranging terms, this immediately implies that
    \begin{IEEEeqnarray}{rCl}
        \int_{\*\Theta \setminus B_{\*\theta_{C}}(r_C)}
        q(\*\theta)d\*\theta
        & < & 
        \dfrac{t+\Delta}{C}
        \nonumber
    \end{IEEEeqnarray}
    Since $C$ was chosen arbitrarily and can be picked arbitrarily large thanks to the coerciveness of $\mathbb{E}_{\mu}\left[\ell(\*\theta, \*x)\right]$, this immediately implies that $\mathcal{S}_t$ is tight.
    Again, by Prokhorov's Theorem this implies that $\mathcal{S}_t$ is compact, which completes the proof.
\end{proof}

\begin{lemma}[Equi-coerciveness]
    If Assumptions \ref{AS:min_exists} and \ref{AS:D} hold, $\{\overline{F}_n\}_{n=1}^{\infty}$ is equi-coercive.
    \label{lemma:equi_coercive_Fbar}
\end{lemma}
%
%
\begin{proof}
    $\{\overline{F}_n\}_{n=1}^{\infty}$ is equi-coercive if and only if there exists a coercive function $\Psi$ for which $\Psi \leq \overline{F}_n$ for all $n$ (see e.g., Proposition 7.7 in \citep[][]{gammaConvergence}).
    Since Assumption \ref{AS:min_exists} ensures that $\mathbb{E}_{q}\left[\mathbb{E}_{\mu}\left[\ell(\*\theta, \*x)\right]\right]$ exists and is finite, it is valid to take $\Psi(q) = \mathbb{E}_{q}\left[
        \mathbb{E}_{\mu}\left[
            \ell(\*\theta, \*x)
        \right]
    \right]$, which clearly yields a lower bound on $\overline{F}_n$ for all $n$. 
    All that remains is to prove that $\Psi$ is coercive, which holds by virtue of Lemma \ref{lemma:coercive_Psi}.
\end{proof}

\begin{lemma}[$\Gamma$-convergence]
    If Assumptions \ref{AS:min_exists}, \ref{AS:suitable} hold,
    $\overline{F}_n(q)$ $\Gamma$-converges to $\mathbb{E}_{q}\left[\mathbb{E}_{\mu}\left[\ell(\*\theta, \*x)\right]\right]$. 
    \label{lemma:gamma_conv_Fbar}
\end{lemma}
\begin{proof}
    %
    %
    Assumption \ref{AS:min_exists} ensures that $\mathbb{E}_{q}\left[\mathbb{E}_{\mu}\left[\ell(\*\theta, \*x)\right]\right]$ either exists and is finite or is infinite.
    In either case, it holds that $\overline{F}_n(q) \leq \overline{F}_{n-1}(q)$ so that $\overline{F}_n$ is a decreasing sequence of functions. 
    Moreover, it is clear that pointwise (i.e. for fixed $q$), $\overline{F}_n(q) \to \mathbb{E}_{q}\left[\mathbb{E}_{\mu}\left[\ell(\*\theta, \*x)\right]\right]$. This holds trivially if $\overline{F}_n = \infty$ for all $n$ and for the finite-valued case provided that $D(q||\pi) < \infty$, which holds by Assumption \ref{AS:suitable}.
    Taken together, this implies that $\overline{F}_n$ $\Gamma$-converges to the lower-semicontinuous envelope of its pointwise limit by Proposition 5.7 in \citep{gammaConvergence}.
    Now since $\mathbb{E}_{q}\left[\mathbb{E}_{\mu}\left[\ell(\*\theta, \*x)\right]\right]$ is itself lower semi-continuous, it is its own lower-semicontinuous envelope. 
    This completes the proof.
\end{proof}
With the necessary lemmas now in hand, one can state the convergence of $\bar{q}_n$ -- the solution sequence corresponding to the auxiliary objective --
to a dirac delta at the $\mu$-population optimal value of $\*\theta$. 
\begin{corollary}[Consistency]
    If Assumptions \ref{AS:min_exists}, \ref{AS:dirac}, \ref{AS:D} and \ref{AS:suitable} hold,
    $\bar{q}_n \overset{\mathcal{D}}{\to}\delta_{\*\theta^{\ast}}$, i.e. the minimizers of $\overline{F}_n$ weakly converge to a point mass at $\*\theta^{\ast}$ as $n\to\infty$. 
    Moreover, $F_n(\overline{q}_n) \to \mathbb{E}_{\mu}\left[\ell(\*\theta^{\ast}, \*x)\right]$ as $n\to\infty$.
    \label{corollary:consistency_Fbar}
\end{corollary}
\begin{proof}
    This is a simple application of Corollary 7.24 in \citep{gammaConvergence}. 
    By Lemmas \ref{lemma:equi_coercive_Fbar} and \ref{lemma:gamma_conv_Fbar}, $\overline{F}_n$ is both equi-coercive and $\Gamma$-convergent to $\mathbb{E}_{q}\left[\mathbb{E}_{\mu}\left[\ell(\*\theta, \*x)\right]\right]$. 
    Further, by Assumptions \ref{AS:min_exists} and \ref{AS:dirac}, the infimum of $\mathbb{E}_{q}\left[\mathbb{E}_{\mu}\left[\ell(\*\theta, \*x)\right]\right]$ over $\mathcal{Q}^{\*\theta}$ is unique and given by $\delta_{\*\theta^{\ast}}$.
    This completes the proof.
\end{proof}

\subsection{$\varepsilon$-minimizers for $\overline{F}_n$}
\label{sec:varepsilon-minimizers}

This section proves that the minimizers $q_n$ are $\varepsilon_n$-minimizers of $\overline{F}_n$. 
%
%
Specifically, Lemma \ref{lemma:Fubini} guarantees that $q_n$ corresponds to a $\mu$-almost surely finite objective value for all $n$. Similarly, it can directly be used to show that the sequence $\varepsilon_n$ consists only of $\mu$-almost surely finite-valued terms.
Lemma \ref{lemma:eps_minimizers} applies Lemma \ref{lemma:Fubini} to derive an explicit form for $\varepsilon_n$. Crucially, this form does \textit{not} depend on $q_n$, but on $\overline{q}_n$. This will turn out to substantially ease remaining proofs: Unlike $q_n$ which depends on $x_{1:n}$ and thus is a  random measure, $\overline{q}_n$ is non-stochastic. 
%


\begin{lemma}
    If Assumptions
    \ref{AS:min_exists}
    and \ref{AS:finite_solution_exists} hold, then it also holds that
    \begin{itemize}
        \item[(i)]
        $\mathbb{E}_{\pi}\left[
            \mathbb{E}_{\mu}
            \left[
            |\ell(\*\theta, \*x)|
            \right]
        \right] < \infty$
        and
        $\mathbb{E}_{\overline{q}_n}\left[
            \mathbb{E}_{\mu}
            \left[
            |\ell(\*\theta, \*x)|
            \right]
        \right] < \infty$;
        \item[(ii)]
            $\mathbb{E}_{\pi}\left[
            \mathbb{E}_{\mu}
            \left[
            \ell(\*\theta, \*x)
            \right]
        \right] = \mathbb{E}_{\mu}\left[
            \mathbb{E}_{\pi}
            \left[
            \ell(\*\theta, \*x)
            \right]
        \right]$ and 
            $\mathbb{E}_{\overline{q}_n}\left[
            \mathbb{E}_{\mu}
            \left[
            \ell(\*\theta, \*x)
            \right]
        \right] = \mathbb{E}_{\mu}\left[
            \mathbb{E}_{\overline{q}_n}
            \left[
            \ell(\*\theta, \*x)
            \right]
        \right]$;
        \item[(iii)]
        $\mathbb{E}_{\pi} \left[ \ell(\*\theta, x_i)\right] < \infty$
        and
        $\mathbb{E}_{\overline{q}_n} \left[ \ell(\*\theta, x_i)\right] < \infty$ $\mu$-almost surely.
    \end{itemize}
    for any $n\in\mathbb{N}$. 
    \label{lemma:Fubini}
\end{lemma}
\begin{proof}
    To unify notation and avoid proving the same things separately for $\pi$ and $\overline{q}_n$, write $\overline{q}_0 = \pi$.
    
    \textbf{(i)}
    First, observe that by Assumption \ref{AS:finite_solution_exists} and by the definition of the objective $\overline{F}_n$, it holds that
    \begin{IEEEeqnarray}{rCl}
        \infty >
        \mathbb{E}_{\overline{q}_0}\left[
            \mathbb{E}_{\mu}
            \left[
            \ell(\*\theta, \*x)
            \right]\right]
        \geq 
        \dots 
        \geq
        \mathbb{E}_{\overline{q}_n}\left[
            \mathbb{E}_{\mu}
            \left[
            \ell(\*\theta, \*x)
        \right]\right]
        \geq
        \mathbb{E}_{\overline{q}_{n+1}}\left[
            \mathbb{E}_{\mu}
            \left[
            \ell(\*\theta, \*x)
        \right]\right]
        \geq
        \dots
        \label{eq:lemma3_pi_bounds}
    \end{IEEEeqnarray}
    %
    It remains to show that this also holds if one takes the absolute value of the loss.
    Denoting by ${1}_A$ the indicator function and by $\mu$ the probability measure as defined via Assumption \ref{AS:min_exists}. Recall that $\mu$ is the probability measure of a random variable $\*x$ on the measure space $(\Omega, \mathcal{F})$ so that $\*x:\Omega \to \mathcal{X}$.
    With this in mind, it holds that
    \begin{IEEEeqnarray}{rCl}
        &&\mathbb{E}_{\overline{q}_n}\left[
            \mathbb{E}_{\mu}
            \left[
            \ell(\*\theta, \*x)
            \right]
        \right] 
        \nonumber \\
        & = &
        \int_{\*\Theta}\int_{\mathcal{X}}\ell(\*\theta, \*x)d\mu(\*x) \overline{q}_n(\*\theta)
        \nonumber \\
        & = &
        \int_{\*\Theta}\int_{\mathcal{X}}\ell(\*\theta, \*x)
        \cdot 1_{\{\ell(\*\theta, \*x) \leq 0\}}(\*x)
        d\mu(\*x) \overline{q}_n(\*\theta)
        + 
        \int_{\*\Theta}\int_{\mathcal{X}}\ell(\*\theta, \*x)\cdot 1_{\{\ell(\*\theta, \*x) > 0\}}(\*x)
        d\mu(\*x) \overline{q}_n(\*\theta)
        \nonumber
    \end{IEEEeqnarray}
    Because $\overline{q}_n \geq 0$, this also immediately implies that one can compute the absolute expectation via  
    \begin{IEEEeqnarray}{rCl}
        &&\mathbb{E}_{\overline{q}_n}\left[
            \mathbb{E}_{\mu}
            \left[
            |\ell(\*\theta, \*x)|
            \right]
        \right] 
        \nonumber \\
        & = &
        -\int_{\*\Theta}\int_{\mathcal{X}}\ell(\*\theta, \*x)
        \cdot 1_{\{\ell(\*\theta, \*x) \leq 0\}}(\*x)
        d\mu(\*x) \overline{q}_n(\*\theta)
        + 
        \int_{\*\Theta}\int_{\mathcal{X}}\ell(\*\theta, \*x)\cdot 1_{\{\ell(\*\theta, \*x) > 0\}}(\*x)
        d\mu(\*x) \overline{q}_n(\*\theta),
        \quad\quad
        \label{eq:proof_lemma_4}
    \end{IEEEeqnarray}
    which will be finite if both integrals by themselves are finite.
    As it turns out, this is indeed the case: by virtue of Assumption \ref{AS:finite_solution_exists} and eq. \eqref{eq:lemma3_pi_bounds}, 
    $\mathbb{E}_{\overline{q}_n}\left[
            \mathbb{E}_{\mu}
            \left[
            \ell(\*\theta, \*x)
            \right]
        \right] < \infty$.
    Moreover, by Assumption \ref{AS:min_exists}, $\mathbb{E}_{\mu}
            \left[
            \ell(\*\theta, \*x)
            \right]$ 
    is bounded below by $\mathbb{E}_{\mu}
            \left[
            \ell(\*\theta^{\ast}, \*x)
            \right] = C < \infty$ so that
    it also holds that
    \begin{IEEEeqnarray}{rCl}
        \int_{\*\Theta}\int_{\mathcal{X}}\ell(\*\theta, \*x)
        \cdot 1_{\{\ell(\*\theta, \*x) \leq 0\}}(\*x)
        d\mu(\*x) \overline{q}_n(\*\theta)
        & \leq &
        \min\{0, C\} < \infty.
        \nonumber
    \end{IEEEeqnarray}
    Thus, the only remaining term in eq. \eqref{eq:proof_lemma_4} must also be finite and (i) follows.

    \textbf{(ii)} By virtue of \textbf{(i)}, one may apply the Fubini-Tonelli Theorem to conclude that $\mathbb{E}_{\overline{q}_n}\left[
            \mathbb{E}_{\mu}
            \left[
            \ell(\*\theta, \*x)
            \right]
        \right]  = 
    \mathbb{E}_{\mu}\left[
            \mathbb{E}_{\overline{q}_n}
            \left[
            \ell(\*\theta, \*x)
            \right]
        \right] <\infty$.

    \textbf{(iii)} 
    By definition of the expectation, it is clear that $\mathbb{E}_{\mu}\left[
            \mathbb{E}_{\overline{q}_n}
            \left[
            \ell(\*\theta, \*x)
            \right]
        \right] <\infty$ if and only if $\mathbb{P}_{\mu}\left(
            \mathbb{E}_{\overline{q}_n}
            \left[
            \ell(\*\theta, \*x)
            \right] = \infty
        \right) = 0$, or equivalently if $\mathbb{P}_{\mu}\left(
            \mathbb{E}_{\overline{q}_n}
            \left[
            \ell(\*\theta, \*x)
            \right] < \infty
        \right) = 1$. 
        In other words, $\mathbb{E}_{\overline{q}_n}
            \left[
            \ell(\*\theta, x_i)
            \right] < \infty$ holds $\mu$-almost surely.
\end{proof}


\begin{lemma}[$\varepsilon_n$-minimizers]
    If Assumptions \ref{AS:min_exists}, \ref{AS:dirac} and \ref{AS:finite_solution_exists} hold, then the sequence $\{q_n\}_{n=1}^{\infty}$ produces finite valued objectives, i.e. $F_n(q_n)<\infty$ or $F_n(q_n,q_n^z) < \infty$ in the presence of latent variables. Moreover, $q_n$ is a $\varepsilon_n$-solution of $\overline{F}_n$ for all $n$, i.e. 
    \begin{IEEEeqnarray}{rCl}
        \overline{F}_n(q_n) 
        & \leq &
        \inf_{q\in\mathcal{Q}^{\*\theta}}\overline{F}_n(q) + \varepsilon_n
        \nonumber
    \end{IEEEeqnarray}
    for a sequence $\{\varepsilon_n\}_{n=1}^{\infty}$ (with 
    $\varepsilon_n < \infty$ $\mu$-almost surely) given by 
    \begin{IEEEeqnarray}{rCl}
        \varepsilon_n & = &
        2\left| 
           \mathbb{E}_{\overline{q}_n}\left[
                \frac{1}{n}\sum_{i=1}^n\ell(\*\theta, x_i) - \mathbb{E}_{\mu}\left[\ell(\*\theta, \*x)\right]
           \right]
        \right|
        \nonumber
    \end{IEEEeqnarray}
    \label{lemma:eps_minimizers}
\end{lemma}
\begin{proof}
    \textbf{$\mu$-almost surely finite-valued objectives:}
    For the case without latent variables, this immediately follows by Lemma \ref{lemma:Fubini}. Recall that the Lemma implies that $\mathbb{E}_{\pi}\left[ \ell(\*\theta, x_i)\right] < \infty$ $\mu$-almost surely, which means that $F_n(q_n) \leq F_n(\pi) < \infty$, $\mu$-almost surely for all $n$.
    The case with latent variables is similarly simple. In particular, observe that if one knew the unobserved components $z_{s(i)}$, then one could directly apply Lemma \ref{lemma:Fubini} again in the same vein as before. 
    Note further that (i) Assumption \ref{AS:dirac} implies that $q_n^z = \delta_{z_{S(n)}}$ is a feasible choice because it is a limit point of the set $\mathcal{Q}_n^z$ and (ii) Assumption \ref{AS:min_exists} guarantees that $\mathbb{E}_{\delta_{z_{S(n)}}}\left[\ell(\*\theta, x_i, \*z_{s(i)})\right] = \ell(\*\theta, x_i, z_{s(i)})$ is finite. 
    From this, one may conclude that $F_n(q_n, q_n^z) \leq F_n(\pi, \delta_{z_{S(n)}}) < \infty$, $\mu$-almost surely for all $n$.
    Here, the last inequality again follows from Lemma \ref{lemma:Fubini}.
    
    \textbf{Finite-valued $\varepsilon_n$:}
    First, define the difference between $F_n(q)$ and $\overline{F}_n(q)$ without latent variables as
    \begin{IEEEeqnarray}{rCl}
        e_n(q) 
        & = &
        \int_{\*\Theta}q(\*\theta)\left[\frac{1}{n}\sum_{i=1}^n\ell(\*\theta, x_i) - \mathbb{E}_{\mu}\left[\ell(\*\theta, \*x)\right]\right]d\*\theta.
        \nonumber
    \end{IEEEeqnarray}
    It is clear that $\varepsilon_n$ is finite-valued if and only if $e(\overline{q}_n)$ is.
    Now, notice that
    \begin{IEEEeqnarray}{rCl}
    e_n(\overline{q}_n)
    & \leq &
        \frac{1}{n}\sum_{i=1}^n
        \underbrace{
        \int_{\*\Theta}\mathbb{E}_{\overline{q}_n}\left[\ell(\*\theta, x_i)\right]
        }_{<\infty, \text{ Lemma \ref{lemma:Fubini}}}
        - 
        \underbrace{
        \mathbb{E}_{\mu}\left[\ell(\*\theta^{\ast}, \*x)\right]
        }_{<\infty, \text{ Assumption \ref{AS:min_exists}}}.
        \nonumber
    \end{IEEEeqnarray}
    Next, define the difference between $F_n(q,p)$ and $\overline{F}_n(q)$ in the presence of latent variables as
    \begin{IEEEeqnarray}{rCl}
        e_n(q, p) 
        & = &
        \int_{\*\Theta}q(\*\theta)\left[
            \frac{1}{n}\sum_{i=1}^n\mathbb{E}_p\left[ \ell(\*\theta, x_i^o, \*z_{s(i)}) \right] - \mathbb{E}_{\mu}\left[\ell(\*\theta, \*x)\right]
        \right]d\*\theta.
        \nonumber
    \end{IEEEeqnarray}
    In this case, is clear that $\varepsilon_n$ is finite-valued if and only if $e_n(q_n, \delta_{z_{S(n)}})$ is. 
    Using the same arguments as in the case without latent variables and recalling that $x_i = (x_i^o, z_{s(i)})$,
    \begin{IEEEeqnarray}{rCl}
    e_n(\overline{q}_n, \delta_{z_{S(n)}})
    & \leq &
        \frac{1}{n}\sum_{i=1}^n
        \underbrace{
        \int_{\*\Theta}\mathbb{E}_{\overline{q}_n}\left[\ell(\*\theta, x_i)\right]
        }_{<\infty, \text{ Lemma \ref{lemma:Fubini}}}
        - 
        \underbrace{
        \mathbb{E}_{\mu}\left[\ell(\*\theta^{\ast}, \*x)\right]
        }_{<\infty, \text{ Assumption \ref{AS:min_exists}}}.
        \nonumber
    \end{IEEEeqnarray}
    
    \textbf{$\varepsilon_n$-solution:}
    Considering again first the case without latent variables, note that 
    \begin{IEEEeqnarray}{rCCCCCCCl}
        \overline{F}_n(q_n) + e_n(q_n) 
        & = &
        F_n(q_n) 
        & = &
        \inf_{q \in \mathcal{Q}^{\*\theta}}\left[\overline{F}_n(q) + e_n(q)\right] 
        & \leq & \overline{F}_n(\overline{q}_n) + e_n(\overline{q}_n),
        \label{eq:lemma_eps_min_IE1}
    \end{IEEEeqnarray}
   $\mu$-almost surely. Further, by definition of $\overline{q}_n$ as the minimizer of $\overline{F}_n$, it also holds that $\overline{F}_n(q_n) \geq  \overline{F}_n(\overline{q}_n)$ $\mu$-almost surely, so that one may conclude that
    \begin{IEEEeqnarray}{rCCCl}
        0 & \leq \overline{F}_n(q_n) - \overline{F}_n(\overline{q}_n)
        & \leq & e_n(\overline{q}_n) - e_n(q_n),
        \label{eq:lemma_eps_minimizers_eq1}
    \end{IEEEeqnarray}
    $\mu$-almost surely, from which it clearly follows that
    \begin{IEEEeqnarray}{rCl}
        e_n(q_n) & \leq & e_n(\overline{q}_n), \label{eq:lemma_eps_minimizer_eq2}
    \end{IEEEeqnarray}
    $\mu$-almost surely. This allows to conclude that the first result indeed follows: $F_n(q_n)<\infty$. Moreover, it also implies that 
    \begin{IEEEeqnarray}{rCCCCCCl}
        0 &\leq &
        e_n(\overline{q}_n) - e_n(q_n)
        & \leq & 
        |e_n(\overline{q}_n)| + |e_n({q}_n)| 
        & \leq & 2|e_n(\overline{q}_n)|,
        \nonumber
    \end{IEEEeqnarray}
    $\mu$-almost surely. With this last result in hand, one can now define the sequence
    \begin{IEEEeqnarray}{rCl}
        \varepsilon_n &= & 2|e_n(\overline{q}_n)| <\infty
        \nonumber
    \end{IEEEeqnarray}
    which together with eq. \eqref{eq:lemma_eps_minimizers_eq1} yields that indeed, 
    \begin{IEEEeqnarray}{rCl}
        \overline{F}_n(q_n) \leq \inf_{q\in\mathcal{Q}^{\*\theta}}\overline{F}_n(q) + \varepsilon_n.
        \nonumber
    \end{IEEEeqnarray}
    
    The case with latent variables is similar: It holds that $\mu$-almost surely,
         \begin{IEEEeqnarray}{rCCCCCCCl}
        \overline{F}_n(q_n) + e_n(q_n, q_n^z) 
        & = &
        F_n(q_n, q_n^z) 
        & \leq &
        F_n(\overline{q}_n, \delta_{z_{S(n)}})
        & = &
        \overline{F}_n(\overline{q}_n) + e_n(\overline{q}_n, \delta_{z_{S(n)}}).
        \label{eq:lemma_eps_min_IE1_latent}
    \end{IEEEeqnarray}
    As before, $\delta_{z_{S(n)}}$ is the dirac measure at the {true} (though in practice unknown) realizations $z_{S(n)}$ of $\*z_{S(n)}$.
    The inequality in eq. \eqref{eq:lemma_eps_min_IE1_latent} holds by definition of $q_n, q_n^z$ and because $\mathcal{Q}^{z}_n$ admits a dirac delta as a limit point due to Assumption \ref{AS:dirac}. 
    %
    %
    From here on out, the proof is virtually identical to the one of Lemma \ref{lemma:eps_minimizers}: Combining the fact that $\overline{F}_n(\overline{q}_n) \leq \overline{F}_n(q_n)$ with eq. \eqref{eq:lemma_eps_min_IE1_latent} yields that
    \begin{IEEEeqnarray}{rCl}
        0 
        \leq
        \overline{F}_n({q}_n)
        -
        \overline{F}_n(\overline{q}_n)
        \leq
        e_n(\overline{q}_n, \delta_{z_{S(n)}})
        - 
        e_n({q}_n, q^z_n),
        \nonumber
    \end{IEEEeqnarray}
    $\mu$-almost surely.
    An immediate consequence is that $\mu$-almost surely,
    \begin{IEEEeqnarray}{rCl}
        e_n({q}_n, q^z_n) \leq e_n(\overline{q}_n, \delta_{z_{S(n)}}),
        \nonumber
    \end{IEEEeqnarray}
    so that one can follow the same steps as in the proof for Lemma \ref{lemma:eps_minimizers} to conclude that  
    \begin{IEEEeqnarray}{rCl}
        \overline{F}_n(q_n) 
        & \leq &
        \inf_{q\in\mathcal{Q}^{\*\theta}}\overline{F}_n(q) + \varepsilon_n
        \nonumber
    \end{IEEEeqnarray}
    as claimed.
\end{proof}

\subsection{Convergence of $\varepsilon_n$-minimizers to zero}
\label{sec:varepsilon-convergence}

While the last section established that $q_n$ are $\varepsilon_n$-minimizers of $\overline{q}_n$, this insight does not suffice to prove consistency unless one can show that $\varepsilon_n$ goes to zero $\mu$-almost surely (in the $\mu$-probability limit). 
Moreover, recall that $\varepsilon_n$ going to zero $\mu$-almost surely (in the $\mu$-probability limit) is equivalent to Assumption \ref{AS:varepsilon_convergence} holding true.
Deriving conditions ensuring that Assumption \ref{AS:varepsilon_convergence} holds true is the technically most demanding part of establishing \GVI consistency. This section addresses various strategies to achieve this goal.
%

The generic strategy for this is as follows:
Invoking Corollary \ref{corollary:consistency_Fbar} and specifically that $F_n(\overline{q}_n) \to \mathbb{E}_{\mu}\left[\ell(\*\theta^{\ast}, \*x)\right]$ as $n\to\infty$, it clearly holds that 
\begin{IEEEeqnarray}{rCl}
    && \lim_{n\to\infty}\int_{\*\Theta}
    \overline{q}_n(\*\theta)
    \left[
        \frac{1}{n}\sum_{i=1}^n\ell(\*\theta, x_i) - 
        \mathbb{E}_{\mu}\left[\ell(\*\theta, \*x)\right]
    \right]d\*\theta
    \nonumber \\
    & = & 
    \lim_{n\to\infty}\int_{\*\Theta}
    \overline{q}_n(\*\theta)
    \left[
        \frac{1}{n}\sum_{i=1}^n\ell(\*\theta, x_i) 
    \right]d\*\theta
    -
    \mathbb{E}_{\mu}\left[ \ell(\*\theta^{\ast}, \*x)\right].
    \nonumber
\end{IEEEeqnarray}
It immediately follows that $\varepsilon_n$ goes to zero $\mu$-almost surely (in the $\mu$-probability limit) if and only if
%
%
\begin{IEEEeqnarray}{rCCCCCl}
    \lim_{n\to\infty}\int_{\*\Theta}\overline{q}_n(\*\theta)\left[\frac{1}{n}\sum_{i=1}^n\ell(\*\theta, x_i)\right]d\*\theta
    & = &
    \int_{\*\Theta}\lim_{n\to\infty}\left\{\overline{q}_n(\*\theta)\left[\frac{1}{n}\sum_{i=1}^n\ell(\*\theta, x_i)\right]\right\}d\*\theta 
    & = & 
    \mathbb{E}_{\mu}\left[ \ell(\*\theta^{\ast}, \*x)\right]. \quad\quad 
    \label{eq:vareps_convergence_criterion}
\end{IEEEeqnarray}
Again, this convergence is meant to occur $\mu$-almost surely (in the $\mu$-probability limit). 
While the second equality follows trivially, the first equality corresponds to pulling the limit operation under the integral and is in general very difficult to establish.
%
In this sense, the main challenge of the current section is to derive conditions under which the first equality holds. That is to say conditions sufficient to prove the convergence of an indefinite integral over a generally unbounded random function $\overline{q}_n(\*\theta)\frac{1}{n}\sum_{i=1}^n\ell(\*\theta, x_i)$ to the integral over its deterministic pointwise limit.

Establishing this result directly can be done in two different ways. 
Firstly, via boundedness conditions for $\ell$. While these boundedness conditions can take different flavours, they all work by ensuring that one may apply a Dominated Convergence Theorem.
Secondly, via a Law of Large Numbers (\LLN) on the triangular array $\{\{\mathbb{E}_{\overline{q}_n}\left[\mathbb{E}_{\mu}\left[\ell(\*\theta, \*x)\right]\right]\}_{i=1}^n\}_{n=1}^{\infty}$. While this is in general more difficult, it is straightforward if the observations correspond to independent and identically distributed random variables.
The remainder of this section explains and elaborates on these two strategies.

%
%
%

\subsubsection{Convergence via boundedness of $\ell$}

While the boundedness conditions presented next can sometimes be restrictive, they have two clear advantages: Firstly, they are often easy to verify in practice. Secondly, they do not depend on strong assumptions on the Data Generating Mechanism for $\*x_{1}, \*x_{2}, \dots$ such as independence. In particular, boundedness conditions allow for the random variables generating the observation sequence to be strongly dependent on one another.
%
\begin{lemma}
    If Assumption \ref{AS:min_exists},  \ref{AS:dirac}, \ref{AS:D}, \ref{AS:suitable}, \ref{AS:finite_solution_exists} hold and (a) or (b) is true where
    \begin{itemize}
        \item[(a)] $\ell(\*\theta, \*x) \leq h(\*x)$ $\mu$-almost surely for all $\*\theta \in \*\Theta$ so that $\mathbb{E}_{\mu}[h(\*x)]<\infty$;
        \item[(b)] $\frac{1}{n}\sum_{i=1}^n\ell(\*\theta, x_i)$ is asymptotically uniformly $\mu$-almost surely bounded, i.e.
        \begin{IEEEeqnarray}{rCl}
            \mathbb{P}_{\mu}\left(             \lim_{n\to\infty}\sup_{\*\theta \in \*\Theta}
                        \Big|\frac{1}{n}\sum_{i=1}^n\ell(\*\theta, x_i) - \mathbb{E}_{\mu}\left[
                            \ell(\*\theta, \*x)
                        \right]\Big| \leq h(\*x)
                    \right) & = & 1,
            \nonumber
        \end{IEEEeqnarray}
        where $\mathbb{E}_{\mu}[h(\*x)]<\infty$.
        Notice that for $h(\*x) = 0$, this is equivalent to requiring a strong uniform law of large numbers to hold for $\frac{1}{n}\sum_{i=1}^n\ell(\*\theta, x_i)$ over $\*\Theta$;
    \end{itemize}
    then, it also follows that $\mu$-almost surely,
    \begin{IEEEeqnarray}{rCl}
        \lim_{n\to\infty}\int_{\*\Theta}\overline{q}_n(\*\theta)\left[\frac{1}{n}\sum_{i=1}^n\ell(\*\theta, x_i)\right]d\*\theta
        & = &
        \mathbb{E}_{\mu}\left[ \ell(\*\theta^{\ast}, \*x)\right] 
        \nonumber
    \end{IEEEeqnarray}
    which is equivalent to saying that Assumption \ref{AS:varepsilon_convergence}(B) holds.
    \label{lemma:trivial_varepsilon_convergence}
\end{lemma}
\begin{proof}
    %
    \textbf{(a):}
    Notice that $\mathbb{E}_{\overline{q}_n}\left[\mathbb{E}_{\mu}[h(\*x)]\right] = \mathbb{E}_{\mu}[h(\*x)]$,
    so the dominated convergence theorem implies that $\lim_{n\to\infty}\int_{\*\Theta}\overline{q}_n(\*\theta)\left[\frac{1}{n}\sum_{i=1}^n\ell(\*\theta, x_i)\right]d\*\theta
        =
        \mathbb{E}_{\mu}\left[ \ell(\*\theta^{\ast}, \*x)\right]$, $\mu$-almost surely.
        
    \textbf{(b):}
    By assumption, for all $n \geq N$ for some $N<\infty$, it holds that $\sup_{\*\theta \in \*\Theta}\Big|\frac{1}{n}\sum_{i=1}^n\ell(\*\theta, x_i) - \mathbb{E}_{\mu}\left[ \ell(\*\theta, \*x)\right]\Big| \leq M + h(\*x)$ $\mu$-almost surely for some constant $M<\infty$. Using again the dominated convergence theorem, one concludes that $\lim_{n\to\infty}\int_{\*\Theta}\overline{q}_n(\*\theta)|\frac{1}{n}\sum_{i=1}^n\ell(\*\theta, x_i) - \mathbb{E}_{\mu}\left[\ell(\*\theta, \*x)\right]|d\*\theta
        = 0$, $\mu$-almost surely.
    By Scheff\'e's Lemma, this implies that $\lim_{n\to\infty}\int_{\*\Theta}\overline{q}_n(\*\theta)\left[\frac{1}{n}\sum_{i=1}^n\ell(\*\theta, x_i)\right]d\*\theta
        =
        \mathbb{E}_{\mu}\left[ \ell(\*\theta^{\ast}, \*x)\right]$, $\mu$-almost surely, which completes the proof.
\end{proof}
\begin{remark}
    Notice that certain special cases for condition (a) are  easily checked. 
    For instance, any absolutely bounded loss function $\ell(\*\theta, \*x) \leq M$ almost surely for all $\*\theta \in \*\Theta$ (such as Huber-losses) trivially satisfies condition (ii).
    Condition (b) is more generally applicable, but harder to check. For $h(\*x) = 0$ for instance, it amounts to establishing a strong uniform law of large numbers. This typically involves proving (stochastic) equicontinuity conditions for $\ell$, see e.g. Chapter 21 in \citet{StochasticLimitTheory}.
\end{remark}

Once the boundedness conditions can be verified, consistency immediately follows. 
The following Corollary states this formally using a number of easily verifiable conditions that are special cases of the ones in in Lemma \ref{lemma:trivial_varepsilon_convergence}.
This way, one obtains a non-exhaustive list of useful example conditions that are easy to check and additionally apply to the case of dependent (e.g., Time Series) data.

\begin{corollary}
    Suppose that for a given \GVI problem $P(\ell, D, \mathcal{Q})$, Assumptions \ref{AS:min_exists}, \ref{AS:dirac}, \ref{AS:D}, \ref{AS:suitable} and \ref{AS:finite_solution_exists} hold. Suppose that additionally, one of the following holds: 
    \begin{itemize}
        \item[(a)] $\ell(\*\theta, \*x) \leq M$ for some constant $M$, $\mu$-almost surely;
        \item[(b)] $\ell(\*\theta, \*x)$ is jointly continuous in $\*\theta$ and $\*x$, $\mu$-almost surely and both $\mathcal{X}$ and $\*\Theta$ are compact;
        \item[(c)] $\sum_{i=1}^n\ell(\*\theta, x_i)$ converges to $\mathbb{E}_{\mu}\left[\ell(\*\theta, \*x)\right]$ both $\mu$-almost surely and uniformly over $\*\Theta$.
    \end{itemize}
    Then, $q_n \overset{D}{\to} \delta_{\*\theta^{\ast}}$ $\mu$-almost surely, i.e. the solution sequence of $P(\ell, D, \mathcal{Q})$ is strongly consistent.
    \label{corollary:summary_consistency_dependence}
\end{corollary}
\begin{proof}
    Conditions (a) and (b) are just special cases of the first condition  in Lemma \ref{lemma:trivial_varepsilon_convergence}. This is trivial for (a). For (b), joint continuity and compactness imply that one can apply the Extreme Value Theorem, which immediately entails that $\ell(\*\theta, \*x) \leq M$ holds for some constant $M$. Thus, (b) is really only a special case of (a). Lastly, (c) is a direct restatement of the second condition in Lemma \ref{lemma:trivial_varepsilon_convergence} with $h = 0$.
    To conclude the proof, it now suffices if one can show that the conditions of Lemma \ref{lemma:trivial_varepsilon_convergence} imply consistency. 
    This trivially holds as they guarantee that $\varepsilon_n$ goes to zero $\mu$-almost surely, which in turn yields consistency by Corollary 7.24 in \citep{gammaConvergence}.
\end{proof}

\subsubsection{Convergence via Laws of Large Numbers (\LLN{}s)}
\label{sec:varepsilon-convergence-LLN}

The conditions for convergence presented thus far rely on some form of  boundedness, but this is not needed. 
In fact, a more general way in which one can establish the integral's convergence is through a  \LLN (for triangular arrays). 
The current section elaborates on this theme. First, Lemma \ref{lemma:LLN_implies_varepsilon_convergence} shows that proving showing a \LLN{} to hold suffices for consistency. Lemma \ref{lemma:restrict_attention_to_A} then delivers an auxiliary result relying on Assumption \ref{AS:A_exists}. This result is used for proving Lemma \ref{lemma:SLLN_for_iid}, which finally guarantees that Assumption \ref{AS:A_exists} together with the i.i.d. assumption on $\*x_i$ suffices for consistency. 
Lastly, Corollaries \ref{corollary:summary_consistency} and \ref{corollary:example_weak_convergence} give minor variations on the same theme. They show that in principle, \LLN{}s can also hold and used for proving consistency if $\*x_i$ is not i.i.d.

\begin{lemma}
    Suppose a Law of Large Numbers (\LLN) holds for the triangular array
    $\{\{
            Z_i^{(n)}
        \}_{i=1}^n\}_{n=1}^{\infty}$  where
    $Z_i^{(n)} = \mathbb{E}_{\overline{q}_n}\left[ \ell(\*\theta, x_i) \right]$. I.e., for $S_n = \frac{1}{n}\sum_{i=1}^n Z_i^{(n)} = \mathbb{E}_{\overline{q}_n}\left[\frac{1}{n}\sum_{i=1}^n\ell(\*\theta, x_i) \right]$, 
    \begin{itemize}
        \item[(A)] $S_n\overset{\mu-\mathcal{P}}{\to} \mathbb{E}_{\mu}[\ell(\*\theta^{\ast}, \*x)]$ or 
        \item[(B)] $S_n \overset{\mu-a.s.}{\to} \mathbb{E}_{\mu}[\ell(\*\theta^{\ast}, \*x)]$.
    \end{itemize}
    If additionally, Assumptions \ref{AS:min_exists}, \ref{AS:dirac}, \ref{AS:D} and \ref{AS:finite_solution_exists} hold,  Assumption \ref{AS:varepsilon_convergence}(A) or \ref{AS:varepsilon_convergence}(B) follows with the same mode of convergence that applies to $S_n$.
    \label{lemma:LLN_implies_varepsilon_convergence}
\end{lemma}
\begin{proof}
    First, notice that Lemma \ref{lemma:Fubini} implies that $\mathbb{E}_{\overline{q}_n}\left[
            \mathbb{E}_{\mu}
            \left[
            \ell(\*\theta, \*x)
            \right]
        \right] = \mathbb{E}_{\mu}\left[
            Z_i^{(n)}
        \right] < \infty$.
   Hence, one can construct the new triangular array $Y_i^{(n)}$ with partial sums $S_n'$, where
   \begin{IEEEeqnarray}{rCl}
       Y_i^{(n)} = Z_i^{(n)} - \mathbb{E}_{\mu}\left[
            Z_i^{(n)}
        \right]; \quad\quad
        S_n' = \frac{1}{n}\sum_{i=1}^nY_i^{(n)}.
        \nonumber
   \end{IEEEeqnarray}
   Noting that Corollary \ref{corollary:consistency_Fbar} implies that $\mathbb{E}_{\mu}\left[
            Z_i^{(n)}
        \right] \to \mathbb{E}_{\mu}\left[\ell(\*\theta^{\ast}, \*x) \right]$ as $n\to\infty$,
        it immediately follows that $S_n' \overset{\mu-\mathcal{P}}{\toL} 0$ if (A) holds and $S_n' \overset{\mu-a.s.}{\toL} 0$ if (B) holds.
        This is equivalent to saying that Assumption \ref{AS:varepsilon_convergence}(A) or Assumption \ref{AS:varepsilon_convergence}(B) holds and completes the proof.
\end{proof}

The next lemma will be useful in conjunction with Lemma \ref{lemma:LLN_implies_varepsilon_convergence}. 
In words, it says that if Assumption \ref{AS:A_exists} holds with a compact set $A \subset \*\Theta$, then it suffices to prove the integral in eq. \eqref{eq:vareps_convergence_criterion} converges over $A$.
This is useful because convergence of definite integrals over compact sets is much easier to establish than proving convergence of indefinite integrals.
This also explains why Assumption \ref{AS:A_exists} is attractive: It is not only interpretable, but also makes proving consistency substantially easier.

\begin{lemma}
    If Assumption \ref{AS:min_exists},  \ref{AS:dirac}, \ref{AS:D}, \ref{AS:suitable}, \ref{AS:finite_solution_exists} hold and Assumption \ref{AS:A_exists} holds with a compact set $A\subset \*\Theta$, 
    \begin{IEEEeqnarray}{rCl}
        \lim_{n\to\infty}\int_{\*\Theta}\overline{q}_n(\*\theta)\left[\frac{1}{n}\sum_{i=1}^n\ell(\*\theta, x_i)\right]d\*\theta
        & = &
        \lim_{n\to\infty}\int_{\*\Theta}1_A(\*\theta)\overline{q}_n(\*\theta)\left[\frac{1}{n}\sum_{i=1}^n\ell(\*\theta, x_i)\right]d\*\theta
        \quad \text{ $\mu$-almost surely.}
        \nonumber
    \end{IEEEeqnarray}
    Consequently, Assumption \ref{AS:varepsilon_convergence}(B) holds (i.e., $\varepsilon_n \overset{\mu-a.s.}{\toL} 0$) if and only if 
    \begin{IEEEeqnarray}{rCl}
        \lim_{n\to\infty}\int_{\*\Theta}1_A(\*\theta)\overline{q}_n(\*\theta)\left[\sum_{i=1}^n\ell(\*\theta, x_i)\right]d\*\theta
        & = &
        \mathbb{E}_{\mu}\left[ \ell(\*\theta^{\ast}, \*x)\right]
        \quad \text{ $\mu$-almost surely.}
    \nonumber
    \end{IEEEeqnarray}
    \label{lemma:restrict_attention_to_A}
\end{lemma}
\begin{proof}
    This follows by the Dominated Convergence Theorem: By Assumption \ref{AS:A_exists}, there is $N<\infty$ so that for all $n\geq N$, $\overline{q}_n(\*\theta) \leq \pi(\*\theta)$ for all $\*\theta \in \*\Theta \setminus A$. 
    Moreover, by Assumption \ref{AS:min_exists} and Lemma \ref{lemma:Fubini}, it also holds that
    \begin{IEEEeqnarray}{rCCCCl}
        \lim_{n\to\infty}\int_{\*\Theta\setminus A}\pi(\*\theta)\left[\sum_{i=1}^n\ell(\*\theta, x_i)\right]d\*\theta
        & = &
        \int_{\*\Theta\setminus A}\lim_{n\to\infty}\left\{\pi(\*\theta)\left[\sum_{i=1}^n\ell(\*\theta, x_i)\right]\right\}d\*\theta 
        & = & 0,
        \nonumber
    \end{IEEEeqnarray}
    where the equalities hold $\mu$-almost surely.
    This immediately implies that on $\*\Theta\setminus A$, $\int_{\*\Theta\setminus A}\overline{q}_n(\*\theta)\left[\sum_{i=1}^n\ell(\*\theta, x_i)\right]d\*\theta$ goes to zero $\mu$-almost surely, too.
    Hence, $\varepsilon_n$ goes to zero $\mu$-almost surely and the result follows.
\end{proof}

Finally, we obtain the first generically applicable result. While it requires Assumption \ref{AS:A_exists} to hold with some compact set $A \subset \*\Theta$, as Section \ref{sec:assumptions} in the main paper explains, this assumption is benign and holds for most cases of interest.

\begin{lemma}
    If $x_i \overset{iid}{\sim} \*x$ and  Assumptions \ref{AS:min_exists},   \ref{AS:dirac}, \ref{AS:D},  \ref{AS:finite_solution_exists},  \ref{AS:A_exists} hold, the triangular array $\{\{
            Z_i^{(n)}
        \}_{i=1}^n\}_{n=1}^{\infty}$ satisfies a strong law of large numbers.
    \label{lemma:SLLN_for_iid}
\end{lemma} 
\begin{proof}
    %
    %
    Under independence, $\{\{Z_i^{(n)}\}_{i=1}^n\}_{n=1}^{\infty}$ is a triangular array with independent and indentically distributed columns in the sense of \citet{TALLN_columns}. 
    In fact, Theorem 2 in the same paper  follows.
    To apply this Theorem, three conditions need to be satisfied. First, note that
    \begin{IEEEeqnarray}{rCl}
        \lim_{n\to\infty}Z_{i}^{(n)}
        = \lim_{n\to\infty}\int_{A}\overline{q}_n(\*\theta)\ell(\*\theta, \*x) d\*\theta
        {=}
         \int_{A}\lim_{n\to\infty}\{\overline{q}_n(\*\theta)\}\ell(\*\theta, \*x) d\*\theta
         =
        \ell(\*\theta^{\ast}, \*x),
        \quad \text{ $\mu$-almost surely}
        \nonumber
    \end{IEEEeqnarray}
    %
    Here, the first equality follows by Lemma \ref{lemma:restrict_attention_to_A} and the second due to the dominated convergence theorem: because the limit does not depend on $i$, so that by compactness of $A$ and the Extreme Value Theorem, there exists $\widetilde{\*\theta}$ such that $\ell(\*\theta, \*x) \leq \ell(\widetilde{\*\theta}, \*x) < \infty$ for any $\*\theta \in A$ and for any realization of $\*x$, $\mu$-almost surely. Notice that the Extreme Value Theorem may be applied because of the continuity requirement on $\ell$ in Assumption \ref{AS:min_exists}. The last inequality then follows by weak convergence of $\overline{q}_n$ to $\delta_{\*\theta^{\ast}}(\*\theta)$.
    Thus, the first requirement of Theorem 2 \citep{TALLN_columns} is satisfied.
    Second, note that
    \begin{IEEEeqnarray}{rCl}
        \lim_{n\to\infty}\mathbb{E}_{\mu}\left[ 
            Z_i^{(n)}
        \right]
        & = &
        \lim_{n\to\infty}
            \int_{A}{\overline{q}_n}(\*\theta)\mathbb{E}_{\mu}\left[\ell(\*\theta, \*x)\right]d\*\theta.
        =
        \mathbb{E}_{\mu}\left[\ell(\*\theta^{\ast}, \*x)\right]
        \nonumber
    \end{IEEEeqnarray}
    where the first equation holds by combining Lemma \ref{lemma:Fubini} with Lemma \ref{lemma:restrict_attention_to_A} and the second equality is a consequence of Corollary \ref{corollary:consistency_Fbar}.
    Third, notice that 
    \begin{IEEEeqnarray}{rCl}
        \mathbb{E}_{\mu}\left[ 
            \sup_{n\geq 1}|Z_{1}^{(n)} - 
            \ell(\*\theta^{\ast}, \*x)|
        \right] 
        & \leq & 
        \mathbb{E}_{\mu}\left[ \sup_{n\geq 1}|Z_{1}^{(n)}|\right]
        + \mathbb{E}_{\mu}\left[ |\ell(\*\theta^{\ast}, \*x)|\right]
        \nonumber \\
        & = &
        \mathbb{E}_{\mu}\left[ \mathbb{E}_{\pi}\left[\ell(\*\theta, \*x) \right] \right]
        + \mathbb{E}_{\mu}\left[ |\ell(\*\theta^{\ast}, \*x)|\right]
        < \infty
    \end{IEEEeqnarray}
    as a direct consequence of Assumption \ref{AS:min_exists} Lemma \ref{lemma:Fubini} and eq. \eqref{eq:lemma3_pi_bounds}.
\end{proof}

The next Corollary collects this last Lemma together with the implications of Lemma \ref{lemma:trivial_varepsilon_convergence}. Thus, it is a collection of the main findings for the case of $\mu$-almost sure convergence of the \GVI posteriors.

\begin{corollary}
    Suppose that for a given \GVI problem $P(\ell, D, \mathcal{Q})$, Assumptions \ref{AS:min_exists}, \ref{AS:dirac}, \ref{AS:D}, \ref{AS:suitable} and \ref{AS:finite_solution_exists} hold. Suppose that additionally, one of the following holds: 
    \begin{itemize}
        \item[(a)] $x_i \overset{iid}{\sim} \*x$ and Assumption \ref{AS:A_exists} holds;
        \item[(b)] A strong law of large numbers holds for the triangular array $\left\{\left\{\mathbb{E}_{\overline{q}_n}\left[ \ell(\*\theta, x_i) \right]\right\}_{i=1}^n\right\}_{n=1}^{\infty}$;
    \end{itemize}
    Then, $q_n \overset{D}{\to} \delta_{\*\theta^{\ast}}$ $\mu$-almost surely, i.e. the parameter posteriors produced by $P(\ell, D, \mathcal{Q})$ are consistent.
    \label{corollary:summary_consistency}
\end{corollary}
\begin{proof}
    If one can show that conditions (a) and (b) imply that $\varepsilon_n$ goes to zero $\mu$-almost surely, then Corollary 7.24 in \citep{gammaConvergence} immediately implies that the result holds.
    Clearly, (a) ensures that $\varepsilon_n$ goes to zero $\mu$-a.s. as $n\to\infty$ by application of Lemma \ref{lemma:SLLN_for_iid}. Similarly, (b) ensures it via Lemma \ref{lemma:LLN_implies_varepsilon_convergence}. 
\end{proof}

Since all results thus far established $\mu$-almost sure convergence, the next Corollary is an example of a weaker convergence result. It holds in the $\mu$-probability limit and relies on a weak law of large numbers for $\left\{\left\{\mathbb{E}_{\overline{q}_n}\left[ \ell(\*\theta, x_i) \right]\right\}_{i=1}^n\right\}_{n=1}^{\infty}$.

\begin{corollary}
    Suppose that for a given \GVI problem $P(\ell, D, \mathcal{Q})$, Assumptions \ref{AS:min_exists}, \ref{AS:dirac}, \ref{AS:D}, \ref{AS:suitable} and \ref{AS:finite_solution_exists} hold and that a weak law of large numbers holds for the triangular array $\left\{\left\{\mathbb{E}_{\overline{q}_n}\left[ \ell(\*\theta, x_i) \right]\right\}_{i=1}^n\right\}_{n=1}^{\infty}$.
    Then, $q_n \overset{D}{\to} \delta_{\*\theta^{\ast}}$ in the $\mu$-probability limit, i.e. the solution sequence of $P(\ell, D, \mathcal{Q})$ is weakly consistent.
    \label{corollary:example_weak_convergence}
\end{corollary}
\begin{proof}
    The proof is the same as for (b) in Corollary \ref{corollary:summary_consistency_dependence}. The only difference is that Corollary 7.24 in \citep{gammaConvergence} now holds in the $\mu$-probability limit instead of $\mu$-almost surely.
\end{proof}

\subsection{Deriving consistency: The base family strategy}
\label{sec:base_family_stragey}


This section motivates what the paper will refer to as \textit{base family strategy} for proving the consistency of \GVI posteriors. This strategy is appealing in practice: Suppose one wishes to establish consistency of $q_n$ when using a rich variational family $\mathcal{Q}_2^{\*\theta}$ for which Assumption \ref{AS:A_exists} is hard to verify. As a running example, one may for instance take $\mathcal{Q}_2^{\*\theta}$ to be a mixture of normals or a neural network parameterizing a normal distribution.  
The base family strategy then works as follows:
\begin{itemize}
    \item[(1)] 
    Pick a  base family $\mathcal{Q}_1 = \mathcal{Q}_1^{\*\theta} \times \mathcal{Q}_{n,1}^z$ that is simple to analyze, is known to produce consistent \GVI posteriors and satisfies $\mathcal{Q}_1 \subseteq \mathcal{Q}_2$ for another family $\mathcal{Q}_2 = \mathcal{Q}_2^{\*\theta} \times \mathcal{Q}_{n,2}^z$;
    \item[(2)]
    Apply Lemma \ref{lemma:Q_subset_consistency} to conclude that consistency extends to $\mathcal{Q}_2$.
\end{itemize}
The next Lemma formalizes the validity of this strategy.
%
%


%

\begin{lemma}
    Suppose Assumptions \ref{AS:min_exists}, \ref{AS:dirac}, \ref{AS:D}, \ref{AS:suitable}, \ref{AS:finite_solution_exists} hold for the \GVI problem $P(\ell, D, \mathcal{Q}_1)$ with solution sequence $\{q_{n,1}\}_{n=1}^{\infty}$.
    Further, let that $\mathcal{Q}_1 \subseteq \mathcal{Q}_2$ and $\{q_{n,2}\}_{n=1}^{\infty}$ be the solution sequence of the \GVI problem $P(\ell, D, \mathcal{Q}_2)$.
    %
    %
    If $q_{n,1} \overset{\mathcal{D}}{\to} \delta_{\*\theta^{\ast}}$ $\mu$-almost surely (in the $\mu$-probability limit), then it also holds that $q_{n,2} \overset{\mathcal{D}}{\to} \delta_{\*\theta^{\ast}}$ $\mu$-almost surely (in the $\mu$-probability limit) as $n\to\infty$.
    \label{lemma:Q_subset_consistency}
\end{lemma}
\begin{proof}
    Throughout, the notation is adapted from the proof of Lemma \ref{lemma:eps_minimizers}. 
    The proof is given only for the case of $\mu$-almost sure convergence because the arguments are the exact same for the $\mu$-probability limit case.
    To start, write $\overline{q}_{n,i} = \arginf_{q \in \mathcal{Q}^{\*\theta}_i}\overline{F}_n(q)$ for $i=1,2$. 
    First, consider the case without latent variables.
    Since $\mathcal{Q}_1 \subseteq \mathcal{Q}_2$, it holds that 
    \begin{IEEEeqnarray}{rCCCCCCCl}
        \overline{F}_n(q_{n,2}) + e_n(q_{n,2}) 
        & = & 
        \inf_{q \in \mathcal{Q}_2^{\*\theta}}F_n(q) 
        & \leq &
        \inf_{q \in \mathcal{Q}_1^{\*\theta}}F_n(q)
        & = &
        {F}_n(q_{n,1}).
        \label{eq:lemma_larger_Q_eq1}
    \end{IEEEeqnarray}
    Moreover, as  $q_{n,1} \overset{\mathcal{D}}{\to} \delta_{\*\theta^{\ast}}$ $\mu$-a.s. and since by virtue of Corollary \ref{corollary:consistency_Fbar}, $\overline{q}_{n,1} \overset{\mathcal{D}}{\to} \delta_{\*\theta^{\ast}}$, the sequence $e_n(q_{n,1})$ converges to zero $\mu$-a.s..
    %
    This means that  
    \begin{IEEEeqnarray}{rCl}
        F_n(q_{n,1}) = \overline{F}_n(q_{n,1}) + e_n(q_{n,1}) \leq \overline{F}_n(\overline{q}_{n,1}) + e_n(\overline{q}_{n,1}) =  \inf_{q\in\mathcal{Q}_1^{\*\theta}}\overline{F}_n(q) + e_n(\overline{q}_{n,1}),
        \nonumber
    \end{IEEEeqnarray} 
    where the inequality holds by definition of $\overline{q}_{n,1}$ and $q_{n,1}$.
    Combining this with eq. \eqref{eq:lemma_larger_Q_eq1}, one finds that 
    \begin{IEEEeqnarray}{rCCCCCCCCCCl}
        \overline{F}_n(q_{n,2}) 
        \leq
        \inf_{q\in\mathcal{Q}_1^{\*\theta}}\overline{F}_n(q) + e_n(\overline{q}_{n,1})
        -
        e_n(q_{n,2}) 
        & \leq &
        \inf_{q\in\mathcal{Q}_1}\overline{F}_n(q) + |e_n(\overline{q}_{n,1})|
        +
        |e_n(q_{n,2})|, \quad\quad
        \label{eq:lemma_10_inf}
    \end{IEEEeqnarray}
    where the last inequality is applied to ensure that $\widetilde{\varepsilon}_n = |e_n(\overline{q}_{n,1})|+ |e_n(q_{n,2})|>0$.
    Clearly, $|e_n(\overline{q}_{n,1})| \to 0$ $\mu$-a.s. because ${q}_{n,1} \overset{\mathcal{D}}{\to} \delta_{\*\theta^{\ast}}$ by assumption and $\overline{q}_{n,1} \overset{\mathcal{D}}{\to} \delta_{\*\theta^{\ast}}$ $\mu$-a.s. by virtue of Corollary \ref{corollary:consistency_Fbar}.
    This reduces the problem to showing that $|e_n(q_{n,2})|$ goes to zero as $n\to\infty$.
    
    Re-using notation of the proof for Lemma \ref{lemma:eps_minimizers}, it similarly holds for the latent variable case that the problem reduces to showing that $|e_n(q_{n,2},\delta_{z_{S(n)}})|$ goes to zero as $n\to\infty$. 
    %
    %
    To see this, first observe that Assumption \ref{AS:dirac} holding for $\mathcal{Q}_{1,n}^z $ together with the fact that $\mathcal{Q}_{1,n}^z \subseteq \mathcal{Q}_{2,n}^z$ immediately implies that $\mathcal{Q}_{2,n}^z$ satisfies Assumption \ref{AS:dirac}, too.
    From here on out, the result essentially follows by the same steps as before. 
    Re-using the notation of the proof for Lemma \ref{lemma:eps_minimizers} and taking the tuples $(q_{n,1},q_{n,1}^z)$ and $(q_{n,2},q_{n,2}^z)$ to be the optima of $F_n(q,p)$ relative to $\mathcal{Q}_1$ and $\mathcal{Q}_2$ respectively, it holds that
    \begin{IEEEeqnarray}{rCCCCCCCl}
        \overline{F}_n(q_{n,2}) + e_n(q_{n,2}, q_{n,2}^z) 
        & = & 
        \inf_{q,p \in \mathcal{Q}_2}F_n(q,p) 
        & \leq &
        \inf_{q,p \in \mathcal{Q}_1}F_n(q,p)
        & = &
        {F}_n(q_{n,1}, q_{n,1}^z).
        \nonumber
    \end{IEEEeqnarray}
    By the exact same steps and logic as in the absence of latent variables, it thus follows that
    \begin{IEEEeqnarray}{rCl}
        F_n(q_{n,1}, q_{n,1}^z) = \overline{F}_n(q_{n,1}) + e_n(q_{n,1},q_{n,1}^z) 
        & \leq & 
        \overline{F}_n(\overline{q}_{n,1}) + e_n(\overline{q}_{n,1},\delta_{z_{S(n)}}) 
        \nonumber \\
        &=&
        \inf_{q\in\mathcal{Q}_1}\overline{F}_n(q) + e_n(\overline{q}_{n,1},\delta_{z_{S(n)}}),
        \nonumber
    \end{IEEEeqnarray} 
    as well as
    \begin{IEEEeqnarray}{rCl}
        \overline{F}_n(q_{n,2}) 
        & \leq & 
        \inf_{q\in\mathcal{Q}_1}\overline{F}_n(q) + e_n(\overline{q}_{n,1},\delta_{z_{S(n)}})
        -
        e_n(q_{n,2},\delta_{z_{S(n)}}) 
        \nonumber \\
        & \leq &
        \inf_{q\in\mathcal{Q}_1}\overline{F}_n(q) + |e_n(\overline{q}_{n,1},\delta_{z_{S(n)}})|
        +
        |e_n(q_{n,2},\delta_{z_{S(n)}})|, 
        \nonumber
    \end{IEEEeqnarray}
    Thus, in the latent variable case one has $\widetilde{\varepsilon}_n = |e_n(\overline{q}_{n,1},\delta_{z_{S(n)}})|+ |e_n(q_{n,2},\delta_{z_{S(n)}})|>0$.
    So exactly as before, $|e_n(\overline{q}_{n,1},\delta_{z_{S(n)}})| \to 0$ $\mu$-a.s. because ${q}_{n,1} \overset{\mathcal{D}}{\to} \delta_{\*\theta^{\ast}}$ by assumption and $\overline{q}_{n,1} \overset{\mathcal{D}}{\to} \delta_{\*\theta^{\ast}}$ $\mu$-a.s. by virtue of Corollary \ref{corollary:consistency_Fbar}.
    So as claimed, this reduces the problem to showing that $|e_n(q_{n,2},\delta_{z_{S(n)}})|$ goes to zero as $n\to\infty$.

    From here on out, the proofs are the same, regardless whether latent variables are absent ($x_i = x_i^o$) or present ($x_i = (x_i^o, z_{s(i)})$).
    The proof is stated for the case without latent variables, but can immediately be adapted to the latent variable case by replacing each occurrence of $e(\overline{q}_{n,1})$ by $e(\overline{q}_{n,1}, \delta_{z_{S(n)}})$, $e({q}_{n,1})$ by $e({q}_{n,1}, \delta_{z_{S(n)}})$ and $e({q}_{n,2})$ by $e({q}_{n,2}, \delta_{z_{S(n)}})$.
    \begin{IEEEeqnarray}{lll}
        e_n(\overline{q}_{n,1}) - e_n(q_{n,2}) 
        = && \nonumber \\
        \frac{1}{n}\left(
            \mathbb{E}_{\overline{q}_{n,1}}\left[\sum_{i=1}^n\ell(\*\theta, x_{i})\right] 
            - 
            \mathbb{E}_{q_{n,2}}\left[\sum_{i=1}^n\ell(\*\theta, x_{i})\right]
        \right)
        +
        \left(
            \mathbb{E}_{q_{n,2}}\left[\mathbb{E}_{\mu}\left[\ell(\*\theta, \*x)\right]\right] 
            - 
            \mathbb{E}_{\overline{q}_{n,1}}\left[\mathbb{E}_{\mu}\left[\ell(\*\theta, \*x)\right]
        \right]
        \right).
        \nonumber
    \end{IEEEeqnarray}
    Because $\mathcal{Q}_1^{\*\theta} \subset \mathcal{Q}_2^{\*\theta}$, it follows immediately that 
    \begin{IEEEeqnarray}{rCl}
        0 \leq \mathbb{E}_{q_{n,1}}\left[\sum_{i=1}^n\ell(\*\theta, x_{i})\right] - \mathbb{E}_{q_{n,2}}\left[\sum_{i=1}^n\ell(\*\theta, x_{i})\right] \leq \mathbb{E}_{\overline{q}_{n,1}}\left[\sum_{i=1}^n\ell(\*\theta, x_{i})\right] - \mathbb{E}_{q_{n,2}}\left[\sum_{i=1}^n\ell(\*\theta, x_{i})\right].
        \nonumber
    \end{IEEEeqnarray}
    Moreover, by definition of $\overline{q}_{n,2}$ and $q_{n,2}$, it also follows that $\mathbb{E}_{q_{n,2}}\left[\mathbb{E}_{\mu}\left[\ell(\*\theta, \*x)\right]\right] \geq \mathbb{E}_{\overline{q}_{n,2}}\left[\mathbb{E}_{\mu}\left[\ell(\*\theta, \*x)\right]\right]$.
    These two insights allow us to bound
    \begin{IEEEeqnarray}{rCl}
        e_n(q_{n,2}) 
        & \leq &
        e_n(\overline{q}_{n,1})
        - 
        \left( \mathbb{E}_{\overline{q}_{n,2}}\left[\mathbb{E}_{\mu}\left[\ell(\*\theta, \*x)\right]\right] - 
        \mathbb{E}_{\overline{q}_{n,1}}\left[\mathbb{E}_{\mu}\left[\ell(\*\theta, \*x)\right]\right]
        \right).
        \nonumber
    \end{IEEEeqnarray}
    The last step consists in showing that the second term, i.e. the difference of the two expectations goes to zero ($\mu$-almost surely).
    It suffices to show that they both converge to the same limit ($\mu$-almost surely).
    By assumption, $\mathbb{E}_{\overline{q}_{n,1}}\left[\mathbb{E}_{\mu}\left[\ell(\*\theta, \*x)\right]\right] \overset{\mu-a.s.}{\toL} \mathbb{E}_{\mu}\left[\ell(\*\theta^{\ast}, \*x)\right]$. 
    Moreover, $\mathbb{E}_{\overline{q}_{n,2}}\left[\mathbb{E}_{\mu}\left[\ell(\*\theta, \*x)\right]\right] \to \mathbb{E}_{\mu}\left[\ell(\*\theta^{\ast}, \*x)\right]$ because $\mathcal{Q}_1 \subset \mathcal{Q}_2$ allows application of Corollary \ref{corollary:consistency_Fbar}. In conjunction with the fact that $e_n(\overline{q}_{n,1})$ goes to zero and eq. \eqref{eq:lemma_10_inf}, this immediately implies that $q_{n,2}$ produces a $\widetilde{\varepsilon}_n$-sequence of solutions for $\widetilde{\varepsilon}_n =|e_n(q_{n,1})|+ |e_n(q_{n,2})|$.
    Thus, $\widetilde{\varepsilon}_n$ goes to zero $\mu$-almost surely, which together with Corollary 7.24 in \citep{gammaConvergence} concludes the proof.
\end{proof}



\subsection{Application: \GVI with the mean field normal variational family}

This last section takes the  mean field variational family -- the choice for $\mathcal{Q}^{\*\theta}$ that is arguably most popular in practice -- and applies various findings of the previous sections to it.
Specifically, the two main strategies for proving consistency directly as derived in Appendix \ref{sec:varepsilon-convergence} are applied to obtain Corollaries \ref{corollary:GVI_MFN_Boundedness_I}, \ref{corollary:GVI_MFN_Boundedness_II} and \ref{corollary:GVI_MFN_LLN}.
Next, the \textit{base family strategy} is used to prove Corollary \ref{corollary:base_family_strategy_MFN}, showing that a variety of  variational families that are more sophisticated than the mean field normal family yield consistency as well.

The following Corollary is direct application of the observations made in Appendix \ref{sec:varepsilon-convergence} to the mean field normal family: Under mild regularity conditions, any sequence of \GVI posteriors built on a bounded loss is consistent.
\begin{corollary}[Mean field normal \GVI posterior consistency via boundedness I]
    Suppose Assumptions \ref{AS:min_exists}, \ref{AS:D}, \ref{AS:suitable}, \ref{AS:finite_solution_exists} hold and $\ell$ satisfies conditions (a) or (b) of Lemma \ref{lemma:trivial_varepsilon_convergence}.  
    If $\mathcal{Q}^{\*\theta}$ is the mean field normal family and $\mathcal{Q}^z_n$ satisfies Assumption \ref{AS:dirac},  $q_n$ is strongly consistent.
    That is, $q_n \overset{\mathcal{D}}{\toL} \delta_{\*\theta^{\ast}}$, $\mu$-almost surely.
    \label{corollary:GVI_MFN_Boundedness_I}
\end{corollary}
\begin{proof}
    The mean field normal family satisfies Assumption \ref{AS:dirac}, so this is a simple restatement of Corollary \ref{corollary:summary_consistency_dependence}.
\end{proof}

The reader may note that there are still four Assumptions left that this result depends on.
While Assumptions \ref{AS:min_exists} and \ref{AS:finite_solution_exists} are harmless, they are not verifiable unless one is comfortable making explicit statements about the data generating mechanism.
As for the other Assumptions imposed, there is a wide variety of settings in which they hold trivially.
The next Corollary provides some example situations for which this is the case. Note that this is by no means exhaustive, but merely an illustration of the broad applicability of the results.

\begin{corollary}[Mean field normal \GVI posterior consistency via boundedness II]
    Suppose Assumptions \ref{AS:min_exists} and \ref{AS:finite_solution_exists} hold. Suppose also that $\ell$ satisfies conditions (a) or (b) of Lemma \ref{lemma:trivial_varepsilon_convergence}.
    If $\mathcal{Q}$ is the mean field normal family,
    if $D$ is the Kullback-Leibler divergence, R\'enyi's $\alpha$-divergence with $\alpha \in (0,1)$, the $\alpha$-divergence with $\alpha \in (0,1)$, the $\beta$-divergence with $\beta > 1$, the $\gamma$-divergence with $\gamma > 1$, the Fisher divergence or the total variation distance 
    and if $\pi$ is a normal distribution, then $q_n$ is strongly consistent.
    That is, $q_n \overset{\mathcal{D}}{\toL} \delta_{\*\theta^{\ast}}$, $\mu$-almost surely.
    \label{corollary:GVI_MFN_Boundedness_II}
\end{corollary}
\begin{proof}
    The proof proceeds by showing that this is a special case of Corollary \ref{corollary:GVI_MFN_Boundedness_I} by showing that Assumptions \ref{AS:dirac}, \ref{AS:D} and \ref{AS:suitable} are satisfied for the permissible choices for $D$, $\pi$ and $\mathcal{Q}$.
    Clearly, the mean field normal family satisfies Assumption \ref{AS:dirac}.
    Further, all candidate divergences $D$ listed satisfy the condition of Assumption \ref{AS:D}. 
    Moreover, if $q$ and $\pi$ are both normals, then $D(q\|\pi)$ can be written down in closed form (see Appendix \ref{appendix:experiments} or the Appendix in \citet{GVI} for some examples). Moreover, these closed forms are finite for any parameterization of $q$ and $\pi$, which is to say that $D(q\|\pi)<\infty$. Hence, Assumption \ref{AS:suitable} is satisfied, too.
\end{proof}

The conditions of the last Corollary on $\mathcal{Q}$ and $D$ are extremely mild, making clear that \GVI consistency will virtually always hold if $\ell$ satisfies sufficiently strong boundedness conditions.
More often than not however, the required conditions on $\ell$ do not hold or are next to impossible to verify.
In this situation, one can rely on the approach outlined in Appendix \ref{sec:varepsilon-convergence-LLN}. For this approach, one needs to establish a Law of Large Numbers (\LLN). In practice, this can be done by combining Assumption \ref{AS:A_exists} with the supposition that the random variables $\*x_i$ generating the observations $x_i$ are independent and identically distributed. 
This idea summarized in Corollary \ref{corollary:summary_consistency} and applied to the mean field normal family next.
The most important part is the establishment of sufficiently weak conditions under which Assumption \ref{AS:A_exists} holds for the mean field normal family.

\begin{lemma}
    Suppose Assumptions \ref{AS:min_exists}, \ref{AS:dirac}, \ref{AS:D}, \ref{AS:suitable} and \ref{AS:finite_solution_exists} hold with $x_i \overset{iid}{\sim} \*x_1$.
    %
    Suppose also that $\mathcal{Q}^{\*\theta}$ is the mean field variational family.
    If $\pi$ has tails that are monotonically decreasing at most as fast as Gaussian tails outside of some compact set $A'$, then Assumption \ref{AS:A_exists} holds.
    \label{lemma:A_exists_for_MFN}
\end{lemma}
\begin{proof}
    First, fix notation and write
    $\mathcal{Q}^{\*\theta} = \{\prod_{d=1}^D\mathcal{N}(\*\theta_d|\*\mu_d, \*\sigma_d): \*\mu \in \mathbb{R}^D, \*\sigma \in \mathbb{R}_{+}^D, \*\theta \in \*\Theta \}$.
    Further, note that $\*\Theta$ is a normed space and that
    as $\overline{q}_n$ is indexed by the variational parameters $\*\kappa_n$, weak convergence to a point mass can only hold if $\*\kappa_n \to \*\kappa^{\ast}$.
    In the case of a univariate normal, it is also clear that $\*\kappa_{n,d} = (\*\mu_{n,d}, \*\sigma_{n,d})' \to (\*\theta^{\ast}_d, 0)$ as $n\to\infty$.
    Moreover, $\pi$ has monotonically decreasing tails decaying at most as fast as Gaussian tails outside $A'$.
    Take $d_1 = \max_{x,y \in A'}||x-y||_2$, $d_2 = \min_{x \in A'}||\*\theta^{\ast} - x||_2$ and set $x^{\ast} = \argmin_{x \in A'}||\*\theta^{\ast} - x||_2$.
    Clearly, $d_1$ has the interpretation of the minimum width of a ball required to enclose all of $A'$, $d_2$ that of the distance between $\*\theta^{\ast}$ and $A'$ an $x^{\ast}$ as the point in $A'$ with minimum distance to $\*\theta^{\ast}$.
    With this, one can construct a set $A$ containing $\*\theta^{\ast}$ in its interior while also, $\pi$ has monotonically decreasing tails decaying at most as fast as Gaussian tails outside of $A$.
    Specifically, take $A = \{\*\theta: ||\*\theta - x^{\ast}||_2 \leq 2\cdot \max\{d_1,d_2\} \}$.
    Further, set $\pi^{\min} = \min_{\*\theta \in \partial A}\pi(\*\theta)$, where $\partial A$ is the boundary of $A$. 
    Because of the tail requirement on $\pi$, if one can find $N$ such that for all $n\geq N$,
    (i) $\int_{A}\overline{q}_n(\*\theta)d\*\theta \geq \int_{A}\pi(\*\theta)d\*\theta$ and (ii) $\overline{q}_n^{\max} = \max_{\*\theta \in \partial A}\overline{q}_n(\*\theta)< \pi^{\min}$, the result follows.
    It is easy to find such $N$ as $N = \max\{N_1, N_2\}$. 
    First, pick $N_1$ large enough so that for all $n\geq N_1$, $\int_{A}\overline{q}_n(\*\theta)d\*\theta \geq \int_{A}\pi(\*\theta)d\*\theta$. Note that this $N_1$ exists since $\overline{q}_n$ weakly converges and $1_A(\*\theta)$ is a bounded and continuous function.
    Second, pick $N_2$ large enough to ensure that the largest value $\overline{q}_n$ takes on the boundary of $A$ is still small enough: I.e., pick $N_2$ so that $\overline{q}_n^{\max} < \pi^{\min}$ for all $n \geq N_2$. Note that this $N_2$ exists because $\overline{q}_n$ converges if and only if $\*\kappa_n$ converges.
    Hence, Assumption \ref{AS:A_exists} is satisfied and the result follows.
    %
\end{proof}
\begin{remark}
    The requirement on the prior $\pi$ is formulated generally here. It holds if $\pi$ is a normal or Student's $t$ distribution. Similarly, it is satisfied if $\pi$ is a mixture of these distributions.
\end{remark}

With this result in hand, one can now replace the boundedness conditions in Corollaries \ref{corollary:GVI_MFN_Boundedness_I} and \ref{corollary:GVI_MFN_Boundedness_II} by assumptions on the data generating mechanism and the prior.

\begin{corollary}[Mean field normal \GVI posterior consistency via \LLN]
    Suppose Assumptions \ref{AS:min_exists}, \ref{AS:dirac}, \ref{AS:D}, \ref{AS:suitable} and \ref{AS:finite_solution_exists} hold with $x_i \overset{iid}{\sim} \*x_1$.
    %
    Suppose also that $\mathcal{Q}^{\*\theta}$ is the mean field variational family.
    If $\pi$ has tails that are monotonically decreasing at most as fast as Gaussian tails outside of some compact set $A'$, then $q_n$ is strongly consistent.
    That is, $q_n \overset{\mathcal{D}}{\toL} \delta_{\*\theta^{\ast}}$, $\mu$-almost surely.
    \label{corollary:GVI_MFN_LLN}
\end{corollary}
\begin{proof}
    This is a simple application of Corollary \ref{corollary:summary_consistency}, which is possible as Lemma \ref{lemma:A_exists_for_MFN} guarantees that Assumption \ref{AS:A_exists} holds.
\end{proof}
\begin{remark}
    As with Corollary \ref{corollary:GVI_MFN_Boundedness_I}, one may wish to replace Assumptions \ref{AS:dirac}, \ref{AS:D} and \ref{AS:suitable} with more interpretable requirements. Once again, numerous mild conditions will suffice to achieve this, for instance the same conditions imposed in Corollary \ref{corollary:GVI_MFN_Boundedness_II}.
\end{remark}

The last part of this section illustrates the usefulness and wide applicability of Lemma \ref{lemma:Q_subset_consistency} in the context of the mean field normal family.
%
If one can establish consistency with respect to the mean field normal family, Lemma \ref{lemma:Q_subset_consistency} can be deployed to immediately extend this result for any superset, which encompasses Neural Networks parameterizing a normal distribution, Gaussian Processes as variational families as well as mixtures of normals.

\begin{corollary}[\GVI posterior consistency via base family strategy with mean field normals]
    Suppose Assumptions \ref{AS:min_exists}, \ref{AS:dirac}, \ref{AS:D}, \ref{AS:suitable}, \ref{AS:finite_solution_exists} and \ref{AS:varepsilon_convergence} hold for the \GVI posteriors based on the mean field normal family, which implies that the posteriors are strongly (weakly) consistent.
    Then, any sequence of \GVI posteriors based on a superset of the mean field normal family is also strongly (weakly) consistent.
    Such supersets include mixtures of normals, Gaussian Processes or a Neural Network parameterizing a normal distribution.
    \label{corollary:base_family_strategy_MFN}
\end{corollary}
\begin{proof}
    Apply Lemma \ref{lemma:Q_subset_consistency}  in conjunction with Lemma \ref{lemma:A_exists_for_MFN}.
\end{proof}

\section{Experiments}
\label{appendix:experiments}

This appendix details derivations and results for the experiments in the main paper.
All experiments are conducted using the black box \GVI procedure outlined by \citet{GVI}.
An implementation is available at \url{https://github.com/JeremiasKnoblauch/GVI_JASA}.

\subsection{Details on \BLR}

The regression coefficients for the \BLR are given by
\begin{IEEEeqnarray}{rCl}
    \beta_{\text{true}} & = & (16.32,  10.15, -12.45,   2.92 ,
         9.21,  -4.20,   5.66,   4.09,
         3.04,   1.25,   7.33 ,  15.03, \nonumber \\
        && -6.65,  13.28,   5.29 ,   7.45,
        -8.37,   4.35,  17.85,  -7.80)^T    
        \nonumber
\end{IEEEeqnarray}

\subsection{Derivation of Uncertainty quantifiers}

A mean field normal variational family $\mathcal{Q}^{\*\theta}$ for $\*\theta = (\*\theta_1, \dots, \*\theta_D)'$ is given by
\begin{IEEEeqnarray}{rCl}
    \mathcal{Q}^{\*\theta} &= & 
    \left\{ 
        q(\*\theta| \*\mu_{q}, \*\sigma_q) = \mathcal{N}(\*\theta| \*\mu_q, \text{diag}^2(\*\sigma_q)): \*\mu_q \in \mathbb{R}^D, \*\sigma_q \in \mathbb{R}_{+}^D
    \right\}, \nonumber
\end{IEEEeqnarray}
where $\text{diag}^2(\*\sigma_q)$ is a $D\times D$ diagonal matrix whose $d$-th entry on the diagonal is $\sigma_{q,d}^2$, i.e. the $d$-th squared entry of $\*\sigma_q$. 
In order for the inference procedure of \citet{GVI} to be applicable, one needs to estimate or compute the uncertainty quantifiers $D(q\|\pi)$. 
For the experiments, the  prior is also a completely factorized normal distribution of form
\begin{IEEEeqnarray}{rCl}
    \pi(\*\theta) & = & \mathcal{N}(\*\theta| \*\mu_{\pi}, \text{diag}^2(\*\sigma_{\pi})) \nonumber
\end{IEEEeqnarray}
In this case (and also if $\pi$ is not factorized), all uncertainty quantifiers used in the experiments are available in closed form.
For the Kullback-Leibler Divergence (\KLD), it is well-known that for two factorized normal distributions $q$ and $\pi$ as defined above, 
\begin{IEEEeqnarray}{rCl}
    \KLD(q\|\pi) & = & 
    \frac{1}{2}\sum_{d=1}^D 
    \left\{
        \frac{\sigma_{q,d}^2}{\sigma_{\pi,d}^2} 
        +
        \frac{(\mu_{d,\pi} - \mu_{j,q})^2}{\sigma_{\pi,j}^2} 
        -
        1
        +
        2\ln\left(\frac{\sigma_{\pi,d}}{\sigma_{q,d}}\right)
    \right\}.
    \nonumber
\end{IEEEeqnarray}
Similarly, R\'enyi's $\alpha$-divergence 
is available in closed form for any $\alpha \in (0,1)$ and any exponential family (see for instance the Appendix in \citet{GVI}). For $q$ and $\pi$ as above, one finds that
\begin{IEEEeqnarray}{C}
    \RAD(q\|\pi) =
    \frac{1}{\alpha-1}
    \sum_{d=1}^D 
    \left\{
        Z_{\alpha,d} - \alpha Z_{q,d} - (1-\alpha) Z_{\pi,d}
    \right\}
    \nonumber
\end{IEEEeqnarray}
In this case, the three quantities to be computed are
\begin{IEEEeqnarray}{C}
    Z_{q,d} = \ln(\sigma_{q,d}^2) + \frac{1}{2}\frac{\mu_{q,d}^2}{\sigma_{q,d}^2}; \quad
    Z_{\pi,d}  =  \ln(\sigma_{\pi,d}^2) + \frac{1}{2}\frac{\mu_{\pi,d}^2}{\sigma_{\pi,d}^2}; \quad
    Z_{\alpha,d} = \ln(\sigma_{\alpha,d}^2) + \frac{1}{2}\mu_{\alpha,d}^2\sigma_{\alpha,d}^2.
    \nonumber
\end{IEEEeqnarray}
In turn, the new quantities used for computing $Z_{\alpha,d}$ are given as
\begin{IEEEeqnarray}{rCl}
    \sigma_{\alpha,d}^2 & = & \left(\frac{\alpha}{\sigma_{q,d}^2} + \frac{1-\alpha}{\sigma_{\pi,d}^2}\right)^{-1}, \quad
    \mu_{\alpha,d} = \alpha \frac{\mu_{q,d}}{\sigma_{q,d}^2} + (1-\alpha) \frac{\mu_{\pi,d}}{\sigma_{\pi,d}^2}
    \nonumber
\end{IEEEeqnarray}
Since the $\alpha$-divergence (\AD) can be rewritten in terms of the \RAD, this also means that the \AD has closed form.
Lastly, the Fisher divergence (\FD) also allows a closed form for this case. It is almost certain that this result has been arrived at before, but a reference could not be located in the literature. 
Thus, the derivations are restated next for completeness' sake.
\begin{proposition}
    For $q \in \mathcal{Q}^{\*\theta}$ and $\pi$ completely factorized normals as above,
    \begin{IEEEeqnarray}{rCl}
        \FD(q\|\pi) & =&
        \sum_{d=1}^D
        \left\{
            C_{1,d}^2 + 2C_{1,d}C_{2,d}\cdot\mu_{q,d} + C_{2,d}^2\cdot(\sigma_{q,d}^2 + \mu_{q,d}^2)
        \right\}
        \nonumber
    \end{IEEEeqnarray}
    where $C_{1,d} = \left(
            \frac{\mu_{q,d}}{\sigma_{q,d}^2} - \frac{\mu_{\pi,d}}{\sigma_{\pi,d}^2} 
        \right)$ 
    and $C_{2,d} = \left( 
            \frac{1}{\sigma_{\pi,d}^2} - \frac{1}{\sigma_{q,d}^2} \right)$.
\end{proposition}
\begin{proof}
    First, simply write out the \FD:
    \begin{IEEEeqnarray}{rCl}
        \FD(q\|\pi) & = &
        \int_{\*\Theta}
        \|
            \nabla_{\*\theta}\log q(\*\theta|\mu_q, \sigma_q) - 
            \nabla_{\*\theta}\log \pi(\*\theta|\mu_q, \sigma_q)
        \|_2^2  \:  
        q(\*\theta|\mu_q, \sigma_q)d\*\theta
        \nonumber.
    \end{IEEEeqnarray}
    As the inside of the squared norm is completely linear, it is easy to see that it suffices to individually compute the $d$-th inside the norm first. With simple derivations, one can show that this term takes the form
    \begin{IEEEeqnarray}{rCl}
        \underbrace{\left(
            \frac{\mu_{q,d}}{\sigma_{q,d}^2} - \frac{\mu_{\pi,d}}{\sigma_{\pi,d}^2} 
        \right)}_{=C_{1,d}} + 
        \*\theta_d
        \underbrace{\left( 
            \frac{1}{\sigma_{\pi,d}^2} - \frac{1}{\sigma_{q,d}^2} \right)}_{=C_{2,d}},
        \nonumber
    \end{IEEEeqnarray}
    which due to the full factorization implies that
    \begin{IEEEeqnarray}{rCl}
        \FD(q\|\pi) & = & 
        \sum_{d=1}^D
        \int_{\*\Theta_d}
        \left(
            C_{1,d} + \*\theta_dC_{2,d}
        \right)^2
        q(\*\theta_d|\mu_{q,d}, \sigma_{q,d}).
        \nonumber
    \end{IEEEeqnarray}
    Because $q(\*\theta_d|\mu_{q,d}, \sigma_{q,d})$ is a normal distribution, its first two moments exist and so
    \begin{IEEEeqnarray}{rCl}
        \int_{\*\Theta_d}
        \left(
            C_{1,d} + \*\theta_dC_{2,d}
        \right)^2
        q(\*\theta_d|\mu_{q,d}, \sigma_{q,d})
        & = &
        C_{1,d}^2 + 2C_{1,d}C_{2,d}\cdot\mu_{q,d} + C_{2,d}^2\cdot(\sigma_{q,d}^2 + \mu_{q,d}^2),
        \nonumber
    \end{IEEEeqnarray}
    from which it follows that
    \begin{IEEEeqnarray}{rCl}
        \FD(q\|\pi) & = & 
        \sum_{d=1}^D
        \left\{
            C_{1,d}^2 + 2C_{1,d}C_{2,d}\cdot\mu_{q,d} + C_{2,d}^2\cdot(\sigma_{q,d}^2 + \mu_{q,d}^2)
        \right\},
    \end{IEEEeqnarray}
    which proves the proposition.
\end{proof}

\end{document}